\documentclass{amsart}
\usepackage{amssymb, enumerate, xspace}
\usepackage[all]{xy}
\SelectTips{cm}{}

\numberwithin{equation}{section}

1

\newenvironment{enumeratei}
{\begin{enumerate}[\upshape (i)]}
{\end{enumerate}}

\newtheorem*{namedtheorem}{\theoremname}
\newcommand{\theoremname}{testing}
\newenvironment{named}[1]{\renewcommand\theoremname{#1}
\begin{namedtheorem}}
{\end{namedtheorem}}

\newtheorem{maintheorem}{Theorem}

\newtheorem{theorem}{Theorem}[section]
\newtheorem{proposition}[theorem]{Proposition}
\newtheorem{proposition-definition}[theorem]
{Proposition-Definition}
\newtheorem{corollary}[theorem]{Corollary}
\newtheorem{lemma}[theorem]{Lemma}

\theoremstyle{definition}
\newtheorem{definition}[theorem]{Definition}
\newtheorem{notation}[theorem]{Notation}

\newtheorem{remark}[theorem]{Remark}

\theoremstyle{remark}

\newcommand\nome{testing}
\newcommand\call[1]{\label{#1}\renewcommand\nome{#1}}
\newcommand\itemref[1]{\item\label{\nome;#1}}
\newcommand\refall[2]{\ref{#1}~(\ref{#1;#2})}
\newcommand\refpart[2]{(\ref{#1;#2})}

\newcommand\cO{\mathcal{O}}

\renewcommand\AA{\mathbb{A}}

\newcommand\PP{\mathbb{P}}

\newcommand\rM{\mathrm{M}}

\newcommand\arr{\ifinner\to\else\longrightarrow\fi}

\newcommand\noqed{\renewcommand\qed{}}

\renewcommand\H{\operatorname{H}}

\newcommand\eqdef{\overset{\mathrm{def}} =}

\newcommand\into{\hookrightarrow}

\newcommand\larrowsim{\overset{\sim}\longrightarrow}

\newcommand\im{\operatorname{im}}

\renewcommand\th{^\text{th}}

\def\displaytimes_#1{\mathrel{\mathop{\times}\limits_{#1}}}

\def\displayotimes_#1{\mathrel{\mathop{\bigotimes}\limits_{#1}}}

\renewcommand\hom{\operatorname{Hom}}

\newcommand\ext{\operatorname{Ext}}

\newcommand\tor{\operatorname{Tor}}

\newcommand\spec{\operatorname{Spec}}

\newcommand\codim{\operatorname{codim}}

\newcommand\generate[1]{\langle #1 \rangle}

\newcommand\rk{\operatorname{rk}}

\newcommand\pr{\operatorname{pr}}

\newdir{ >}{{}*!/-5pt/@{>}}

\newcommand\K[1][*]{\operatorname{K}_{#1}}

\newcommand\spK{\operatorname{\mathbb{K}}}

\newcommand\gm[1][S]{\mathbb{G}_{\mathrm{m},#1}}

\def\chq#1){\ch#1)_{\mathbb{Q}}}

\newcommand\N{\mathrm{N}}

\newcommand\W{\operatorname{W}}

\newcommand\R{\mathrm{R}}

\newcommand\Rtilde{\mathrm{\widetilde R}}

\renewcommand\sp{\operatorname{Sp}}

\newcommand\cone{\operatorname{Cone}}

\renewcommand\star{\operatorname{Star}}

\newcommand\ch[1][*]{\operatorname{A}^{#1}}

\newcommand\cm{Cohen--Macau\-lay\xspace}

\newcommand\geom{_{\operatorname{geom}}}

\newcommand\loc{\operatorname{loc}}

\begin{document}

\title[Higher algebraic K-theory for actions of diagonalizable
groups]{Higher algebraic K-theory for actions\\of diagonalizable
groups}

\author{Gabriele Vezzosi}

\address{Dipartimento di Matematica Applicata\\
Universit\`a di Firenze\\
I-50139 Firenze\\ Italy}
\email{vezzosi@dma.unifi.it}

\author{Angelo Vistoli}

\address{Dipartimento di Matematica\\
Universit\`a di Bologna\\
40126 Bologna\\ Italy}
\email{vistoli@dm.unibo.it}

\begin{abstract}
We study the K-theory of actions of diagonalizable group schemes on
noetherian regular separated algebraic spaces: our main result shows how to
reconstruct the K-theory ring of such an action from the K-theory rings of
the loci where the stabilizers have constant dimension. We apply this
to the calculation of the equivariant K-theory of toric varieties,
and give conditions under which the Merkurjev spectral sequence
degenerates, so that the equivariant K-theory ring determines the ordinary
K-theory ring. We also prove a very refined localization theorem for actions
of this type.
\end{abstract}

\date{September 18, 2004}
\subjclass[2000]{19E08; 14L30}

\thanks{Paper published in Invent. Math. \textbf{153} (2003), 1--44.\\Partially supported by the University of Bologna, funds
for selected research topics.}

\maketitle

\tableofcontents

%%%%%%%%%%%%%%%%%%%%%%%%%%%%%%%%%%

\section*{Introduction}

Fix a basis noetherian separated connected scheme $S$, and let $G$ be a
diagonalizable group scheme of finite type over $S$ (see
\cite[Expos\'e VII]{sga3}); recall that this means that $G$ is the product
of finitely many multiplicative groups
$\gm$ and group schemes $\boldsymbol{\mu}_{n,S}$ of
$n^\mathrm{th}$ roots of 1 for various values of $n$. Suppose that $G$
acts on a separated noetherian regular algebraic space $X$ over $S$.

If $G$ acts on $X$ with finite stabilizers, then \cite{vevi} gives a
decomposition theorem for the equivariant higher K-theory ring $\K(X,G)$;
it says that, after inverting some primes, $\K(X,G)$ is a product of
certain factor rings $\K(X^\sigma,G)_\sigma$ for each subgroup schemes
$\sigma \subseteq G$ with $\sigma \simeq \boldsymbol{\mu}_n$ for some $n$
and $X^\sigma \neq \emptyset$ (the primes to be inverted are precisely the
ones dividing the orders of the $\sigma$). A slightly
weaker version of this theorem was given in \cite{toen}. From this one can
prove analogous formulas assuming that the stabilizers are of constant
dimension (Theorem~\ref{thm:refinedconstant}).

This paper deals with the general case, when the dimensions of the
stabilizers are allowed to jump. In this case one sees already in the
simplest examples that $\K(X,G)$ will not decompose as a product, not even
after tensoring with $\mathbb{Q}$; for example, if $S$ is a field, $G$ is
a torus and $X$ is a representation of $G$, then $\K[0](X,G)$ is the ring
of representations $\R G$, which is a ring of Laurent polynomials over
$\mathbb{Z}$.

However, we show that the ring $\K(X,G)$ has a canonical structure of
fibered product. More precisely, for each integer $s$ we consider the
locus $X_s$ of $X$ where the stabilizers have dimension precisely equal to
$s$; this is a locally closed regular subspace of $X$. For each $s$ consider
the normal bundle $\N_s$ of $X_s$ in $X$, and the subspace $\N_{s,s-1}$
where the stabilizers have dimension precisely $s-1$. There is a pullback
map $\K(X_s,G) \to \K(\N_{s,s-1},G)$; furthermore in
Section~\ref{sec:specializations} we define a specialization homomorphism
$\sp_{X,s}^{s-1} \colon \K(X_{s-1},G) \to \K(\N_{s,s-1},G)$, via
deformation to the normal bundle. Our first main result
(Theorem~\ref{thm:maintheorem}) show that these specializations
homomorphisms are precisely what is needed to reconstruct the equivariant
K-theory of $X$ from the equivariant K-theory of the strata.

\begin{maintheorem}[The theorem of reconstruction
from the strata]\call{main:maintheorem} Let $n$ be the dimension of $G$.
The restriction homomorphisms
    \[
    \K(X,G) \longrightarrow \K(X_s,G)
    \]
induce an isomorphism
    \begin{align*}
    \K(X,G) \simeq {}&\K(X_n,G) \displaytimes_{\K(\N_{n,n-1},G)}
    \K(X_{n-1},G) \displaytimes_{\K(\N_{n-1,n-2},G)}\\
    &\quad\ldots
    \displaytimes_{\K(\N_{2,1},G)}  \K(X_1,G) \displaytimes_{\K(\N_{1,0},G)}
    \K(X_0,G).
    \end{align*}
In other words: the restrictions $\K(X,G) \to \K(X_s,G)$ induce an
injective homomorphism $\K(X,G) \to \prod_s \K(X_s,G)$, and an
element $(\alpha_n, \ldots, \alpha_0)$ of the product $\prod_s \K(X_s,G)$
is in the image of $\K(X,G)$ if and only if the pullback of $\alpha_s \in
\K(X_s,G)$ to $\K(\N_{s,s-1},G)$ coincides with\/
$\sp_{X,s}^{s-1}(\alpha_{s-1}) \in \K(\N_{s,s-1},G)$ for all $s = 1$,
\dots,~$n$.
\end{maintheorem}

This theorem is a powerful tool in studying the K-theory of
diagonalizable group actions. From it one gets
easily a description of the higher equivariant K-theory of regular toric
varieties (Theorem~\ref{thm:describetoric}). This is analogous to the
description of their equivariant Chow ring in \cite[Theorem~5.4]{brion1}.

One can put Theorem~\ref{main:maintheorem} above together with the main
result of \cite{vevi} to give a very refined description of $\K(X,G)$;
this is Theorem~\ref{thm:refineddecomposition}, which can be considered the
ultimate localization theorem for actions of diagonalizable groups.
However, notice that it does not supersede
Theorem~\ref{main:maintheorem}, because Theorem~\ref{main:maintheorem}
holds with integral coefficients, while for the formula of
Theorem~\ref{thm:refineddecomposition} to be correct we have to invert some
primes.

Many results are known for equivariant intersection theory, or for
equivariant cohomology; often one can use our theorem to prove their
K-theoretic analogues. For example, consider the following theorem of Brion,
inspired in turn by results in equivariant cohomology due to Atiyah
(\cite{atiyah}), Bredon (\cite{bredon}), Hsiang (\cite{hsiang}), Chang and
Skjelbred (\cite{ch-sk}), Kirwan (\cite{kirwan}),  Goresky, Kottwitz and
MacPherson (\cite{gkmp}); see also the very useful discussion in
\cite{brion2}.

\begin{named}{Theorem}[{\cite[3.2, 3.3]{brion1}}]\call{thm:brion} Suppose
that $X$ is a smooth projective algebraic variety over an algebraically
closed field with an action of an algebraic torus $G$.

\begin{enumeratei}
\itemref{1} The rational equivariant Chow ring $\ch_G(X)_{\mathbb{Q}}$
is free as a module over $\ch_G(\mathrm{pt})_{\mathbb{Q}}$.

\itemref{2} The restriction homomorphism
    \[
    \ch_G(X)_{\mathbb{Q}} \longrightarrow
    \ch_G(X^G)_{\mathbb{Q}} =
    \ch(X^G)_{\mathbb{Q}} \otimes \ch_G(\mathrm{pt})_{\mathbb{Q}}
    \]
is
injective, and its image is the intersection of all the images of the
restriction homomorphisms $\ch_G(X^T)_{\mathbb{Q}} \to
\ch_G(X^G)_{\mathbb{Q}}$, where $T$ ranges over all the subtori of
codimension 1.
\end{enumeratei}
\end{named}

 From this one gets a very simple description of the rational equivariant
Chow ring when the fixed point locus $X^G$ is zero dimensional, and the
fixed point set $X^T$ is at  most 1-dimensional for any subtorus $T
\subseteq G$ of codimension~1 (\cite[Theorem~3.4]{brion1}).

In this paper we prove a  version of Brion's theorem for algebraic
K-theory. Remarkably, it holds with integral coefficients:
we do not need to tensor with $\mathbb{Q}$. This confirms the authors'
impression that when it comes to torsion, K-theory tends to be better
behaved than cohomology, or intersection theory.

The following is a particular case of
Corollary~\ref{cor:maincorollary}; when $G$ is a torus, it is an an
analogue of part~\refpart{thm:brion}{2} of Brion's theorem.

\begin{maintheorem}\label{main:theorem2}
Suppose that $G$ is a diagonalizable group acting a smooth proper scheme
$X$ over a perfect field; denote by $G_0$ the toral component of $G$,
that is, the largest subtorus contained in $G$.

Then the restriction homomorphism
$\K(X, G) \to  \K(X^{G_0}, G)$ is injective, and its image equals the
intersection of all the images of the restriction homomorphisms $\K(X^T,
G) \to \K(X^{G_0}, G)$ for all the subtori $T \subseteq G$ of
codimension~1.
\end{maintheorem}

 From this one gets a very complete description of $\K(X,G)$ when $G$ is a
torus and $X$ is smooth and proper over an algebraically closed field, in
the ``generic'' situation when $X$ contains only finitely many invariant
points and finitely many invariant curves (Corollary~\ref{cor:generic}).

We also analyze the case of smooth toric varieties in detail in
Section~\ref{sec:toric}.

The analogue of Theorem~\ref{main:theorem2} should hold for the integral
equivariant topological K-theory of a compact differentable manifold with
the action of a compact torus. Some related topological results are
contained in \cite{rokn}.

\subsection*{Description of contents} Section~\ref{sec:notation} contains
the setup that will be used throughout this paper. The K-theory that we use
is the one described in \cite{vevi}: see the discussion in
Subsection~\ref{sub:equiK-theory}.

Section~\ref{sec:preliminary} contains some
preliminary technical results; the most substantial of these is a very
general self-intersection formula, proved following closely the
proof of Thomason of the analogous formula in the non-equivariant case
(\cite[Th\'eo\-r\`eme~3.1]{th1}). Here we also discuss the stratification
by dimensions of stabilizers, which is our basic object of study.

In Section~\ref{sec:specializations} we define various types of
specializations to the normal bundle in equivariant K-theory.
This is easy for $\K[0]$, but for the whole higher K-theory ring we
do not know how to give a definition in general without using the
language of spectra.

Section~\ref{sec:reconstruction} contains the proof of
Theorem~\ref{main:maintheorem}.

Section~\ref{sec:limits} is dedicated to the analysis of the case
when $X$ is complete, or, more generally, admits enough limits
(Definition~\ref{def:admitslimits}). The condition that $\K(X,G)$ be free
as a module over the representation ring $\R G$ is not
adequate when working with integral coefficients: here we analyze a
rather subtle condition on the $\R G$-module $\K(X,G)$ that ensures that
the analogue of Brion's theorem above holds, then we show, using a
Bia\l ynicki-Birula stratification, that this condition is in fact satisfied
when $X$ admits enough limits over a perfect field.

We also apply our machinery to show that the degeneracy of the Merkurjev
spectral sequence in \cite{merk}, that he proves when $X$ is smooth and
projective, in fact happens for torus actions with enough limits.

Section~\ref{sec:toric} is dedicated to the K-theory of smooth toric
varieties. For any smooth toric variety $X$ for a torus $T$, we give two
descriptions of $\K(X,T)$. First of all, we show how
Theorem~\ref{thm:maintheorem} in this case gives a simple description of
it as a subring of a product of representation rings, analogous to the
description of its equivariant Chow ring in \cite[Theorem~5.4]{brion1}.
Furthermore, we give a presentation of $\K(X,T)$ by generators and
relations over the K-theory ring of the base field (Theorem~\ref{thm:SR}),
analogous to the classical Stanley--Reisner presentation for its
equivariant cohomology first obtained in \cite{bdcp}. For $\K[0]$ the result
is essentially stated in \cite{kly}.

In Section~\ref{sec:decomposition} we generalize the result of \cite{vevi}
by giving a formula that holds for all actions of diagonalizable groups on
regular noetherian algebraic spaces, irrespective of the dimensions of the
stabilizers (Theorem~\ref{thm:refineddecomposition}).

\subsection*{Acknowledgments}

The first author is thankful for the hospitality at the University of
Utah and the Universit\'e de Grenoble, where some of the work on this
paper has been done.

During the writing of the first draft of this paper the second author held
a visiting position at the Department of Mathematics of the University of
Utah: he is very grateful for the hospitality. He would like to thank
Jim Carlson, Steve Gersten, Dragan Milicic, Anne and Paul Roberts, Paula and
Domingo Toledo, and, most particularly, Aaron Bertram, Herb Clemens and
their families for making his stay such an enjoyable experience.

Both of us would like to thank Michel Brion for his help, and for many
very helpful conversations; among other things, he suggested the
possibility that we might be able to prove
Corollary~\ref{cor:maincorollary}. Also, his excellent articles
\cite{brion1} and \cite{brion2} were a source of inspiration.

The results of Subsection~\ref{subsec:multSR} were obtained after
stimulating conversations with Bernd Sturmfels, Howard
Thompson, Allen Knutson and Corrado De Concini, that took place while the
second author was a visitor at the Mathematical Sciences Research Institute
in Berkeley. Knutson also pointed out the references \cite{atiyah} and
\cite{bredon} to us. We are grateful to all of them, and to MSRI for the
hospitality.

We also want to thank the referee, who did an unusually careful and fair job.

Finally, we would like to remark how much we owe to the articles of the
late Robert W. Thomason: without his groundbreaking work on K-theory,
equivariant and non, this paper might never have been written. His
premature death has been a great blow to the mathematical community.

\section{Notations and conventions}\label{sec:notation}

Throughout the paper we fix a base scheme $S$, that is
assumed to be connected, separated and noetherian.

We will denote by $G$ a diagonalizable group scheme of finite type over $S$
(see
\cite{sga3}), except when otherwise mentioned. Its groups of characters is
$\widehat G \eqdef \hom_S(G, \gm)$; the contravariant functor from the
category of diagonalizable groups schemes of finite type over $S$ to the
category of finitely generated abelian groups given by $G
\mapsto \widehat G$ is an antiequivalence of categories. The ring of
representations of $G$ is, by definition, $\R G = \mathbb{Z}
\widehat G$, and furthermore $G = \spec \R G \times_{\spec \mathbb{Z}} S$.

We will denote by $G_0$ the \emph{toral part} of $G$, that is, the largest
subtorus of $G$. The group of characters $\widehat G_0$ is the quotient of
$\widehat G$ by its torsion subgroup.

A $G$-space will always be
a regular separated noetherian algebraic space over $S$ over which $G$
acts; sometimes we will talk about a \emph{regular} $G$-space, for
emphasis.

We notice explicitly that if $S' \to S$ is a morphism of schemes, with
$S'$ connected, then every diagonalizable subgroup scheme of $G\times_S
S'$ is obtained by base change from a unique diagonalizable subgroup scheme
of $G$. This will be used as follows: if $p: \spec \Omega \to X$ is a
geometric point, then we will refer to its stabilizer, which is a priori
a subgroup scheme of $G\times_S \spec \Omega$, as a subgroup scheme of $G$.

If $Y \into X$ is a regular embedding,
we denote by $\N_Y X$ the normal
bundle.

\subsection{Equivariant K-theory}\label{sub:equiK-theory}

In this subsection $G$ will be a group scheme over $S$ that is
flat, affine and of finite type. We use the same K-theoretic
setup as in \cite{vevi}, that uses the language of \cite{thtr}. The
following is a slight extension of \cite[Theorem~6.4]{vevi}.

\begin{proposition}\call{prop:K-theory}
Let $G$ be flat affine separated group scheme of finite type over $S$,
acting over a noetherian regular separated scheme $X$ over $S$. Consider
the following complicial bi-Waldhausen categories:
\begin{enumeratei}
\itemref{qc} the category $\W_1(X,G)$ of complexes of quasicoherent
$G$-equivariant  $\mathcal{O}_X$-mod\-ules with bounded coherent
cohomology;

\itemref{c} the category $\W_2(X,G)$ of bounded complexes of coherent
$G$-equivariant
$\mathcal{O}_X$-modules;

\itemref{fqc} the category $\W_3(X,G)$ of complexes of flat quasicoherent
$G$-equivariant  $\mathcal{O}_X$-modules with bounded coherent cohomology,
and

\itemref{fqcba} the category $\W_4(X,G)$ of bounded above
complexes of $G$-equivariant quasi-coherent flat
$\mathcal{O}_X$-Modules with bounded coherent cohomology.
\end{enumeratei}

Then the inclusions
    \[
    \W_2(X,G) \subseteq \W_1(X,G)\quad \text{and}\quad\W_4(X,G)
    \subseteq\W_3(X,G)
    \subseteq \W_1(X,G)
    \]
induce isomorphisms on the corresponding Waldhausen K-theories.
Furthermore the K-theory of any of the categories above coincides with the
Quillen K-theory $\K'(X,G)$ of the category of $G$-equivariant coherent
$\mathcal{O}_X$-modules.
\end{proposition}

\begin{proof}
For the first three categories, and the Quillen K-theory, the statement is
precisely \cite[Theorem~6.4]{vevi}.

Let us check that the inclusion $\W_4(X,G) \subseteq \W_1(X,G)$ induces an
isomorphism in K-theory. By  \cite[Proposition 6.2]{vevi}, which shows
that hypothesis 1.9.5.1 is satisfied, we can apply \cite[Lemma
1.9.5]{thtr}, in the situation where $\mathcal{A}$\/ is the category of
$G$-equivariant quasicoherent $\mathcal{O}_X$-Modules,
$\mathcal{C}$ the category of cohomologically bounded complexes in
$\mathcal{A}$,  $\mathcal{D}$ the category of $G$-equivariant
quasicoherent flat $\mathcal{O}_X$-Modules,
$F:\mathcal{D}  \into \mathcal{A}$ is the natural inclusion. In
particular, any complex in $\W_1(X,G)$ receives a quasi-isomorphism from a
complex in $\W_4(X,G)$. That is,
\cite[1.9.7.1]{thtr}, applied to the inclusion
$\W_4(X,G)\into \W_1(X,G)$, is satisfied; since the other hypothesis
1.9.7.0 of \cite[1.9.7]{thtr} is obviously satisfied, we conclude by
\cite[Theorem 1.9.8]{thtr}.
\end{proof}

We will denote by $\spK(X,G)$ the Waldhausen K-theory spectrum and by
$\K(X,G)$ the Waldhausen K-theory group of any of the categories above. As
observed in \cite[p. 39]{vevi}, it follows from results of Thomason that
$\K(-,G)$ is a covariant functor for proper maps of noetherian regular
separated $G$-algebraic spaces over $S$; furthermore, each $\K(X,G)$ has a
natural structure of a graded ring, and each equivariant morphism $f
\colon X \to Y$ of noetherian regular separated
$G$-algebraic spaces over $S$ induces a pullback $f^* \colon \K(Y,G) \to
\K(X,G)$, making $\K(-,G)$ into a contravariant functor from the category
of noetherian regular separated $G$-algebraic spaces over $S$ to
graded-commutative rings. Furthermore, if $i \colon Y \into X$ is a
closed embedding of  noetherian regular separated $G$-algebraic spaces and
$j \colon X \setminus Y \into X$ is the open embedding, then
$\spK(X\setminus Y, G)$ is the cone of the pushforward map
$i_* \colon \spK(Y,G) \to \spK(X,G)$ (\cite[Theorem~2.7]{th5}),
so there is an exact localization sequence
    \[
    \cdots\longrightarrow \K[n](Y,G) \overset{i_*} \longrightarrow
    \K[n](X,G)
    \overset{j^*} \longrightarrow \K[n](X \setminus Y,G) \overset{\partial}
    \longrightarrow
    \K[n-1](Y,G) \longrightarrow \cdots
    \]

Furthermore, if $\pi\colon E \to X$ is a $G$-equivariant vector bundle,
the pullback
    \[
    \pi^*\colon \K(X,G) \longrightarrow \K(E,G)
    \]
is an isomorphism
(\cite[Theorem~4.1]{th5}).

\section{Preliminary results}\label{sec:preliminary}

\subsection{The self-intersection formula}

Here we generalize Thomason's self-inter\-section formula
(\cite[Th\'eor\`eme~3.1]{th1}) to the equivariant case.

\begin{theorem}[The self-intersection formula]\call{thm:selfintersection}
Suppose that a flat group scheme $G$ separated and of finite type over $S$
acts over a noetherian regular separated algebraic space $X$. Let
$i \colon Z \into X$ be a regular $G$-in\-variant closed subspace of $X$.
Then
   \[
   i^*i_* \colon \spK(Z,G) \longrightarrow
   \spK(Z,G),
   \]
coincides up to homotopy with the cup product
   \[
   \lambda_{-1}(\N_Z^\vee X)\smile(-) \colon \spK(Z,G)
   \longrightarrow \spK(Z,G),
   \]
where $\N_Z^\vee X$ is the conormal sheaf of $Z$ in $X$.

In particular, we have the equality
    \[
    i^*i_*  = \lambda_{-1}(\N_Z^\vee X)\smile(-)
    \colon \K(Z,G)\longrightarrow
    \K(Z,G).
    \]
\end{theorem}

\begin{proof}
The proof follows closely Thomason's proof of \cite[Th\'eor\`eme 3.1]
{th1}, therefore we will only indicate the
changes we need for that proof to adapt to our situation.

Let us denote by $\W'(Z,G)$ the Waldhausen category consisting of
pairs $(E^{\bullet },\lambda :L^{\bullet \bullet }\rightarrow
i_{*}E^{\bullet })$ where $E^{\bullet }$ is a bounded above
complex of $G$-equivariant quasi-coherent flat
$\mathcal{O}_Z$-Modules with bounded
coherent cohomology, $L^{\bullet \bullet }$ is a bicomplex of
$G$-equivariant quasi-coherent flat $\mathcal{O}_X$-Modules with
bounded coherent total cohomology such that $L^{ij}=0$ for $j\leq
0$, any $i$ and
also $L^{ij}=0$ for $i>N$, for some integer $N$, any $i$; finally
$\lambda
:L^{\bullet \bullet }\rightarrow i_{*}E^{\bullet }$ is an exact
augmentation of the bicomplex $L^{\bullet \bullet }$. In
particular, for any $i$, the horizontal complex $L^{i\bullet }$
is a flat resolution of $i_{*}E^i$. The morphisms, cofibrations
and weak equivalences in $\W'(Z,G)$ are as in
\cite[3.3, p. 209]{th1}. Thomason \cite[3.3]{th1} shows that
the forgetful functor $(E^{\bullet },\lambda :L^{\bullet \bullet
}\rightarrow i_{*}E^{\bullet })\mapsto E^{\bullet }$ from
$\W'(Z,G)$ to the category $\W_4(Z,G)$ of bounded above
complexes of $G$-equivariant quasi-coherent flat
$\mathcal{O}_Z$-Modules with bounded coherent cohomology induces
a homotopy equivalence between the associated Waldhausen
$K$-theory spectra. In other words, by Proposition~\ref{prop:K-theory}, we
can (and will) use $\W'(Z,G)$ as a ``model'' for
$\spK(Z,G)$.

With these choices, the morphism of spectra $i^*i_*:\spK(
Z,G) \longrightarrow \spK(Z,G)$ can be
represented by the exact functor $\W'(Z,G)\longrightarrow
\W_1(Z,G)$ which sends $(E^{\bullet },\lambda :L^{\bullet
\bullet }\rightarrow i_{*}E^{\bullet })$ to the total complex of
the bicomplex $i^{*}(L^{\bullet \bullet })$.

The rest of the proof is exactly the same as in \cite[3.3, pp.
210-212]{th1}. One first consider functors $T_k:\W'(Z,G)\to
\W_1(Z,G)$ sending an object $(E^{\bullet },\lambda :L^{\bullet
\bullet }\rightarrow i_{*}E^{\bullet })$ to the total complex of
the bicomplex
    \[
    \xymatrix@-6pt{&{}\vdots\ar[d]&{}\vdots\ar[d]
    &&{}\vdots\ar[d]&{}\vdots\ar[d]\\
    0\ar[r]& {}\im\partial_h^{i,-k-1}\ar[d]\ar[r]& i^{*}L^{i,-k}\ar[d]\ar[r]
    & {}\cdots\ar[r]&
    i^*L^{i,-1}\ar[d]\ar[r] & i^*{L^{i,0}}\ar[d]\ar[r]& 0\\
    0\ar[r]&{} \im\partial_h^{i+1,-k-1}\ar[r]\ar[d]&
    i^{*}L^{i+1,-k}\ar[r]\ar[d]&
    {}\cdots\ar[r]& i^*L^{i+1,-1}\ar[r]\ar[d]
    &{} i^*{L^{i+1,0}}\ar[r]\ar[d]& 0\\
    &{}\vdots&{}\vdots&&{}\vdots&{}\vdots
    }
    \]
which results from truncating all the horizontal
complexes of $i^{*}L^{\bullet \bullet }$ at the $k$-th level.

The functors $T_k$ are zero for $k<0$ and come naturally
equipped with functorial epimorphisms $T_k\twoheadrightarrow
T_{k-1}$ whose kernel $h_k$ has the property that $h_k(E^{\bullet
},\lambda :L^{\bullet \bullet }\rightarrow i_{*}E^{\bullet })$ is
quasi isomorphic to $\Lambda ^k(\mathrm{N}_Z^{\vee }X)\otimes
_{\mathcal{O}_Z}i^{*}E^{\bullet }[k]$ (\cite[3.4.4]{th1},
essentially because each horizontal complex in $L^{i\bullet }$ is
a flat resolution of $i_{*}E^i$. Therefore, by induction on
$k\geq -1$, starting from $T_{-1}=0$, each $T_k$ has values in
$\W_1(Z,G)$ and preserves quasi-isomorphisms. Moreover, the
arguments in \cite[3.4, pp. 211-212]{th1}, show that $T_k$
actually preserves cofibrations and pushouts along cofibrations;
hence each $T_k:\W'(Z,G)\longrightarrow \W_1(Z,G)$ is an
exact functor of Waldhausen categories.

As in \cite[3.4, p. 212]{th1}, the quasi-isomorphism
    \[ h_k(E^{\bullet },\lambda :L^{\bullet \bullet }\rightarrow
    i_{*}E^{\bullet })\simeq \Lambda ^k(\N_Z^{\vee }X)\otimes
    _{\mathcal{O}  _Z}i^{*}E^{\bullet }[k]
    \]
shows that the canonical truncation morphism $i^{*}(L^{\bullet
\bullet })\rightarrow T_d(E^{\bullet },\lambda :L^{\bullet
\bullet }\rightarrow i_{*}E^{\bullet })$, $d$ being the
codimension of $Z$ in $X$, is a quasi-isomorphism, i.e. the
morphism of spectra $i^{*}i_{*}:\spK\left( Z,G\right)
\longrightarrow \spK\left( Z,G\right) $ can also be
represented by the exact functor $T_d:\W'(Z,G)\longrightarrow \W_1(Z,G)$.
Now, the Additivity Theorem  (\cite[1.7.3 and 1.7.4]{thtr}) shows that
the canonical exact sequences of functors $ h_k\into
T_k\twoheadrightarrow T_{k-1}$ yield up-to-homotopy equalities
$T_k=T_{k-1}+h_k$ between the induced map of spectra. And finally,
recalling that a shift $[k]$ induces multiplication by $(-1)^k$ at the
level of spectra, by induction on $k\geq -1$ we get equalities up to
homotopy
    \begin{align*}
    i_{}^{*}i_{*} &= T_d(-)\\
    &= \sum_kh_k(-)\\
    &= \sum_k[\Lambda ^k(\N_Z^{\vee }X)]\otimes
    _{\mathcal{O}_Z}i^{*}(-)[k]\\
    &= \lambda _{-1}(\N_Z^{\vee }X)\smile (-)
    \end{align*}
of morphisms of spectra $\spK(Z,G) \to \spK(Z,G)$ (\cite[p.
212]{th1}).
\end{proof}

\subsection{Stratification by dimensions of
stabilizers}\label{subsec:stratdimension}

Let $G$ be a diagonalizable group scheme of finite type acting on $X$ as
usual. Consider the group scheme $H \to X$ of stabilizers of the action.
Since for a point
$x \in X$ the dimension of the fiber $H_x$ equals its dimension at the
point $\gamma(x)$, where $\gamma\colon X \to H$ is the unit section, it
follows from Chevalley's theorem (\cite[13.1.3]{ega4}) that there is an
open subset $X_{\le s}$ where the fibers of $H$ have dimension at most
$s$. We will use also $X_{<s}$ with a similar meaning. We denote
by $X_s$ the locally closed subset $X_{\le s} \setminus X_{<s}$; we will
think of it as a subspace of $X$ with the reduced scheme structure. Finally,
we call $\N_s$ the normal bundle of $X_s$ in $X$, and $\N^0_s$ the
complement of the 0-section in $\N_s$. Notice that $G$ acts on $\N_s$, so we
may consider the subscheme $(\N_s)_{<s} \subseteq \N_s$.

\begin{proposition}\call{prop:firststratification}
Let $s$ be a nonzero integer.

\begin{enumeratei}

\itemref{3} There exists a finite number of $s$-dimensional subtori
$T_1$, \dots,~$T_r$ in $G$ such that $X_s$ is the disjoint union of the
$X_{\le s}^{T_j}$.

\itemref{1} $X_s$ is a regular locally closed subspace
of $X$.

\itemref{2} $\N_s^0 = (\N_s)_{<s}$.
\end{enumeratei}
\end{proposition}

\begin{proof} To prove part\refpart{prop:firststratification}{3} we may
restrict the action of $G$ to its toral component. By Thomason's generic
slice theorem (\cite[Proposition~4.10]{th2}) there are only finitely many
possible diagonalizable subgroup schemes of $G$ that appear as stabilizers
of a geometric point of $X$. Then we can take the $T_j$ to be the toral
components of the $s$-dimensional stabilizers.

Parts \refpart{prop:firststratification}{1} and
\refpart{prop:firststratification}{2} follow from
\refpart{prop:firststratification}{3} and \cite[Proposition~3.1]{th3}.
\end{proof}

\section{Specializations}\label{sec:specializations}

In this section $G$ will be a flat, affine and separated group scheme of
finite type over $S$, acting on a noetherian regular separated algebraic
space $Y$ over $S$.

\begin{definition} A $G$-invariant morphism $Y \to
\mathbb{P}^1_S$ is
\emph{regular at infinity} if the inverse image $Y_\infty$ of the section
at infinity in
$\mathbb{P}^1_S$ is a regular effective Cartier divisor on $Y$.
\end{definition}

\begin{theorem}\call{thm:specializations}
Let $\pi\colon Y \to \mathbb{P}^1_S$ be a $G$-invariant morphism over
$S$ that is regular at infinity. Denote by $i_\infty \colon Y_\infty
\into Y$ the inclusion of the fiber at infinity, $j_\infty \colon Y
\setminus Y_\infty \into Y$ the inclusion of the complement. Then there
exists a specialization homomorphism of graded groups
    \[
    \sp_Y \colon \K(Y \setminus Y_\infty,G) \longrightarrow \K(Y_\infty,G)
    \]
such that the composition
    \[
    \K(Y,G) \overset{j_\infty^*}\longrightarrow \K(Y \setminus Y_\infty)
    \overset{\sp_Y}\longrightarrow
    \K(Y_\infty,G)
    \]
coincides with $i_\infty^* \colon \K(Y,G) \to \K(Y_\infty,G)$.

Furthermore, if $Y'$ is another noetherian separated regular algebraic
space over $S$ and $f \colon Y' \to Y$ is a $G$-equivariant morphism over
$S$ such that the composition $\pi f \colon Y' \to \mathbb{P}^1_S$ is
regular at infinity, then the diagram
    \[
    \xymatrix{{}\K(Y \setminus Y_\infty,G)\ar[r]^>>>>>{\sp_Y} \ar[d]^{f^*}
    & {}\K(Y_\infty,G)\ar[d]^{f_\infty^*}\\
    {}\K(Y' \setminus Y'_\infty,G) \ar[r]^>>>>>{\sp_{Y'}}&
    {}\K(Y'_\infty,G)}
    \]
commutes; here $f_\infty$ is the restriction of $f$ to $Y'_\infty
  \to Y_\infty$.
\end{theorem}

We refer to this last property as \emph{the compatibility of
specializations}.

\begin{proof}
Let us denote by $\spK(X,G)$ the Quillen K-theory spectrum
associated with the category of coherent equivariant $G$-sheaves on a
noetherian  separated algebraic space $X$. There is a homotopy equivalence
    \[
    \cone\bigl(\spK(Y_\infty,G) \xrightarrow {i_{\infty
    *}} \spK(Y,G)\bigr) \simeq \spK(Y\setminus
    Y_\infty, G).
    \]
The commutative diagram
    \[
    \xymatrix{{}\spK(Y_\infty, G)\ar[r]^{i_{\infty*}} \ar@{=} [d]&
    {}\spK(Y,G) \ar[d]^{i_\infty^*}\\
    {}\spK(Y_\infty, G))\ar[r]^{i_\infty^* i_{\infty*}} &
    {}\spK(Y_\infty, G)}
    \]
induces a morphism of spectra
    \begin{equation}\begin{split}\label{morspectra}
    \spK(Y \setminus Y_\infty, G) {}&\simeq
    \cone\bigl(\spK(Y_\infty,G) \overset {i_{\infty
    *}}\longrightarrow\spK(Y,G)\bigr)\\
    &\longrightarrow
    \cone\bigl(\spK(Y_\infty,G) \overset {i_\infty^*
    i_{\infty*}}\longrightarrow
    \spK(Y_\infty,G)\bigr).
    \end{split}\end{equation}
By the self-intersection formula (Theorem~\ref{thm:selfintersection})
there is a homotopy
     \[
     i_\infty^* i_{\infty*} \simeq
     \lambda_{-1}(\N^\vee_{Y_\infty})\smile(-)\colon
     \spK(Y_\infty,G) \to\spK(Y_\infty,G);
     \]
on the other hand $\lambda_{-1}(N^\vee_{Y_\infty})\smile(-)$ is homotopic
to zero, because $\N^\vee_{Y_\infty}$ is trivial. So we have that
    \[
    \cone\bigl(\spK(Y_\infty,G) \xrightarrow{i_\infty^*
    i_{\infty*}} \spK(Y_\infty,G)\bigr) \simeq
    \spK(Y_\infty,G)[1] \oplus \spK(Y_\infty,G),
    \]
where $(-)[1]$ is the suspension of $(-)$. We define the specialization
morphism of spectra
    \[
    \mathcal{S}_Y \colon \spK(Y \setminus Y_\infty, G)
    \longrightarrow \spK(Y_\infty, G)
    \]
by composing the morphism \eqref{morspectra} with the canonical projection
    \[
    \spK(Y_\infty,G)[1] \oplus \spK(Y_\infty,G)
    \longrightarrow \spK(Y_\infty,G).
    \]
Finally, $\sp_Y$ is defined to be the homomorphism induced by
$\mathcal{S}_Y$ on homotopy groups.

Let us check compatibility; it suffices to show that the diagram of
spectra
    \[
    \xymatrix{{}\spK\bigl(Y'_\infty,G\bigr)\ar[r]^{i'_{\infty*}}
    \ar[d]^{f_\infty^*}
    & {}\spK\bigl(Y',G\bigr)\ar[d]^{f^*}\\
    {}\spK\bigl(Y_\infty,G\bigr)\ar[r]^{i_{\infty*}} &
    {}\spK\bigl(Y,G\bigr)}
    \]
commutes up to homotopy. The essential point is that the diagram of
algebraic spaces
    \[
    \xymatrix{Y'_\infty\ar[r]^{i'_{\infty}}\ar[d]^{f_\infty} & Y'\ar[d]^{f}\\
    Y_\infty\ar[r]^{i_\infty} & Y}
    \]
is Tor-independent, that is, $\tor_i^{\mathcal{O}_Y}(\mathcal{O}_{Y'},
\mathcal{O}_{Y_\infty}) = 0$ for all $i>0$. Write $\W(Y',G)$ and
$\W(Y'_\infty,G)$ for the Waldhausen categories of
$G$-equivariant complexes of quasicoherent $\mathcal{O}_{Y'}$-modules and
$\mathcal{O}_{Y'_\infty}$-modules with bounded coherent cohomology, while
$\W(Y,G)$ and $\W(Y_\infty,G)$ will denote the Waldhausen categories of
complexes of $G$-equivariant quasicoherent $\mathcal{O}_Y$-modules and
$\mathcal{O}_{Y_{\infty}}$ bounded coherent cohomology, that are
respectively degreewise $f^*$-acyclic and degreewise $f_\infty^*$-acyclic.
By the Tor-independence of the diagram above, we have that $i_{\infty *}$
gives a functor $\W(Y_\infty,G) \to \W(Y,G)$, and the diagram
    \[
    \xymatrix{{}\W(Y_\infty, G)\ar[r]^{i_{\infty*}} \ar[d]^{f_\infty^*}&
    {}\W(Y,G) \ar[d]^{f^*}\\
    {}\W(Y'_\infty,G)\ar[r]^{i'_{\infty*}} & {}\W(Y',G)}
    \]
commutes. By \cite[1.5.4]{thtr} this concludes the proof of the theorem.
\end{proof}

\begin{remark}\call{rmrk:stupidspecialization} For the projection
$\mathrm{pr}_2 \colon X \times_S \mathbb{P}^1_S \to \mathbb{P}^1_S$ the
specialization homomorphism
    \[
    \sp_{X \times_S \mathbb{P}^1_S} \colon \K(X \times_S \mathbb{A}^1_S,G)
    \longrightarrow \K(X,G)
    \]
coincides with the pullback $s_0^* \colon \K(X \times_S \mathbb{A}^1_S,G)
    \longrightarrow \K(X,G)$ via the zero-section $s_0 \colon X
\to X \times_S \mathbb{A}^1_S$.

In fact, the pullback $j_\infty^* \colon \K(X \times_S
\mathbb{P}^1_S) \to
\K(X \times_S \mathbb{A}^1_S)$ is surjective, and we have
    \[
    \sp_{X \times_S \mathbb{P}^1_S} \circ j_\infty^*
    = s_0^* \circ j_\infty^* = i_\infty^* \colon  \K(X \times_S
    \mathbb{P}^1_S) \longrightarrow \K(X,G).
    \]
\end{remark}

\begin{remark}
Since the restriction homomorphism $\K[0](Y,G) \to \K[0](Y
\setminus Y_\infty,G)$ is a surjective ring homomorphism, and
its composition with $\sp_Y \colon \K[0](Y \setminus Y_\infty,G) \to
\K[0](Y_\infty,G)$ is a ring homomorphism, it follows that the
specialization in degree 0 $\sp_Y \colon \K[0](Y \setminus Y_\infty,G) \to
\K[0](Y_\infty,G)$ is also a ring homomorphism. This should be true for
the whole specialization homomorphism $\sp_Y \colon \K(Y \setminus
Y_\infty,G) \to \K(Y_\infty,G)$, but this is not obvious from the
construction, and we do not know how to prove it.
\end{remark}

\begin{remark} From the proof of Theorem~\ref{thm:specializations} we see
that one can define a specialization homomorphism $\K(Y \setminus
Z,G) \to
\K(Z,G)$ if $Z$ is a regular effective $G$-invariant divisor on $Y$ whose
normal sheaf is $G$-equivariantly trivial.
\end{remark}

\subsection{Specializations to the normal
bundle}\label{subsec:specializationsnormal}

Let us go back to our standard situation, in which $G$ is a diagonalizable
group scheme of finite type acting on a regular separated noetherian
algebraic space $X$. Fix a nonnegative integer $s$, and consider the
closed immersion $X_s
\into X_{\le s}$; denote by $\N_s$ its normal bundle. Consider the
deformation to the normal cone $\pi\colon \mathrm{M}_s \to \mathbb{P}^1_S$,
the one denoted by $\mathrm{M}^0_{X_s} X_{\le s}$ in
\cite[Chapter~5]{fulton}. The morphism $\pi\colon \mathrm{M}_s
\to \mathbb{P}^1_S$ is flat and $G$-invariant. Furthermore $\pi^{-1}(
\mathbb{A}^1_S) = X_{\le s} \times_S ( \mathbb{A}^1_S)$, while the fiber
at infinity of $\pi$ is
$\N_s$. Consider the restriction $\pi^0 \colon \mathrm{M}_s^0 \to
\mathbb{P}^1_S$ to the open subset $\mathrm{M}_s^0 =
(\mathrm{M}_s)_{<s}$; then $(\pi^0)^{-1}( \mathbb{A}^1_S) = X_{<
s} \times_S \mathbb{A}^1_S$, while the fiber at infinity of $\pi^0$ is
$\N_s^0 = (\N_s)_{<s}$. We define a specialization homomorphism
    \[
    \sp_{X,s} \colon \K(X_{<s},G) \longrightarrow \K(\N_s^0, G)
    \]
by composing the pullback
    \[
    \K(X_{<s},G) \longrightarrow \K(X_{<s}\times_S \mathbb{A}^1_S,G)
=
    \K(\mathrm{M}_s^0 \setminus\N_s^0, G)
    \]
with the specialization homomorphism
    \[
    \sp_{\mathrm{M}_s^0} \colon \K(\mathrm{M}_s^0
    \setminus\N_s^0, G)
    \longrightarrow  \K(\N_s^0, G)
    \]
defined in the previous subsection.

We can also define more refined specializations.

\begin{proposition}\call{prop:restrictregularatinfinity}
Let $Y \to \mathbb{P}^1_S$ be regular at
infinity.
\begin{enumeratei}
\itemref{1} If $H\subseteq G$ is a diagonalizable subgroup scheme, then
the restriction  $Y^H \to \mathbb{P}^1$ is regular at infinity.

\itemref{2} If $s$ is a nonnegative integer, the restriction $Y_s \to
\mathbb{P}^1_S$ is also regular at infinity.
\end{enumeratei}
\end{proposition}

\begin{proof}
Part~\refpart{prop:restrictregularatinfinity}{1} and Proposition
\refall{prop:firststratification}{3} imply
part~\refpart{prop:restrictregularatinfinity}{2}.

To prove \refpart{prop:restrictregularatinfinity}{1}, notice that, by
\cite[Prop.~3.1]{th3}, $Y^H$ is regular, and so is $Y_\infty^H$. Let $f$
be the pullback to $Y$ of a local equation for the section at infinity
of $\mathbb{P}^1_S \to S$, and let $p$ be a point of $Y^H_\infty$. Since
the conormal space to $Y^H$ in $Y$ has no nontrivial $H$-invariants,
clearly the differential of $f$ at $p$ can not lie in this conormal
space, hence $f$ is not zero in any neighborhood of
$p$ in
$Y^H$. This implies that $Y^H_\infty$ is a regular Cartier divisor on
$Y^H$, as claimed.
\end{proof}

If $t$ is an integer with $t < s$, let us set $\N_{s,t} \eqdef
(\N_s)_t$. We have that the restriction $(\mathrm{M}_s)_t \to
\mathbb{P}^1_S$ is still regular at infinity, by
Proposition \refall{prop:restrictregularatinfinity}{2}; so we can also
define a specialization homomorphism
    \[
    \sp_{X,s}^t \colon \K(X_{t},G) \longrightarrow \K\bigl(\N_{s,t},
    G\bigr)
    \]
by composing the pullback
   \[
   \K(X_t,G) \longrightarrow \K(X_t\times_S \mathbb{A}^1_S,G) =
   \K\bigl((\mathrm{M}_s^0)_t \setminus(\N_{s,t}, G\bigr)
   \]
with the specialization homomorphism
    \[
    \sp_{{(\mathrm{M}_s^0})_t} \colon \K\bigl((\mathrm{M}_s^0)_t
    \setminus \N_{s,t}, G\bigr)
    \longrightarrow \K\bigl(\N_{s,t}, G\bigr).
    \]

The specializations above are compatible, in the following sense.

\begin{proposition}\label{prop:compspecializations}
In the situation above, the diagram
    \[
    \xymatrix{{}\K(X_{<s},G)\ar[r]\ar[d]^{\sp_{X,s}}&
    {}\K(X_t,G)\ar[d]^{\sp_{X,s}^t}\\
    {}\K(\N_s^0,G)\ar[r]& {}\K(\N_{s,t},G)},
    \]
where the rows are restriction homomorphisms, commutes.
\end{proposition}

\begin{proof}
This follows immediately from the compatibility of specializations
(Theorem~\ref{thm:specializations}).
\end{proof}

\section{Reconstruction from the strata}\label{sec:reconstruction}

\subsection{K-rigidity}

\begin{definition}
Let $Y$ be a $G$-invariant regular locally closed subspace of $X$. We say
that $Y$ is \emph{K-rigid} inside $X$ if $Y$ is regular and
$\lambda_{-1}(\N_Y^\vee X)$ is not a zero-divisor in the ring $\K(Y,G)$.
\end{definition}

This condition may seem unlikely to ever be verified: in the
non-equivariant case $\lambda_{-1}(\N_Y^\vee X)$ is always a nilpotent
element, since it has rank zero over each component of $X$. However, in the
equivariant case this is not necessary true. Here is the basic criterion
that we will use use to check that a subspace is K-rigid.

\begin{lemma}\call{lem:criterionK-rigid}
Let $Y$ be a  $G$-space, $E$ an equivariant vector bundle on $Y$. Suppose
that there is a subtorus $T$ of $G$ acting trivially on $Y$, such that in
the eigenspace decomposition of $E$ with respect to $T$ the subbundle
corresponding to the trivial character is 0. Then $\lambda_{-1}(E)$ is not
a zero-divisor in $\K(Y,G)$.
\end{lemma}

\begin{proof}
Choose a splitting $G \simeq D\times T$; by \cite[Lemme~5.6]{th4}, we
have
    \[
    \K(Y,G) = \K(Y,D) \otimes \R T = \K(Y,D) \otimes \mathbb{Z}[t_1^{\pm 1},
    \ldots, t_n^{\pm 1}].
    \]
If $E = \bigoplus_{\chi\in \widehat T}E_\chi$ is the eigenspace
decomposition of $E$, we have that $\lambda_{-1}(E)$ corresponds to the
element $\prod_{\chi\in \widehat T}\lambda_{-1}(E_\chi \otimes \chi)$ of
$\K(Y,D) \otimes \R T$, so it enough to show that
$\lambda_{-1}(E_\chi \otimes \chi)$ is not a zero-divisor in $\K(Y,D)
\otimes \R T$. But we can write
    \[
    \lambda_{-1}(E_\chi \otimes \chi) = 1 + r_1\chi + r_2\chi^2+
    \cdots + r_n\chi^n \in \K(Y,D) \otimes \R T,
    \]
where $r_n = (-1)^n[\det E_\chi]$ is a unit in $\K(Y,D)$.

Now we can apply the
following elementary fact: suppose that $A$ is a ring,  $r_1$, \dots,~$r_n$
central elements of $A$ such that $r_n$ is a unit, $\chi \in A[t_1^{\pm 1},
\ldots, t_n^{\pm 1}]$ a monomial different from 1. Then the element $1 +
r_1\chi + r_2\chi^2+ \cdots + r_n\chi^n$ is not a zero-divisor in
$A[t_1^{\pm 1}, \ldots, t_n^{\pm 1}]$.
\end{proof}

The next Proposition is a K-theoretic variant of
\cite[Proposition~{3.2}]{brion2}.

\begin{proposition}\call{prop:K-rigid->}
Let $Y$ be a closed K-rigid subspace of $X$, and set $U = X \setminus
Y$. Call $i \colon Y \into X$ and $j \colon U \into X$ the inclusions.

\begin{enumeratei}

\itemref{1} The sequence
    \[
    0 \longrightarrow \K(Y,G) \overset {i_*} \longrightarrow
    \K(X,G)\overset{j^*} \longrightarrow \K(U,G) \longrightarrow 0
    \]
is exact.

\itemref{2} The two restriction maps
    \[
    i^* \colon
    \K(X,G) \longrightarrow \K(Y,G)\quad \mbox{and} \quad j^* \colon \K(X,G)
    \longrightarrow \K(U,G)
    \]
induce a ring  isomorphism
    \[
    (i^*,j^*) \colon \K(X,G) \overset{\sim}\longrightarrow \K(Y,G)
    \displaytimes_{\K(Y,G)/(\lambda_{-1}(\N_Y^\vee X))}
    \K(U,G)
    \]
where $\N_Y^\vee X$ is the conormal bundle of $Y$ in $X$, the homomorphism
    \[
    \K(Y,G) \longrightarrow\K(Y,G)/(\lambda_{-1}(\N_Y X))
    \]
is the projection, while the homomorphism
    \[
    \K(U,G) \simeq \K(X,G)/i_*\K(Y,G) \longrightarrow
    \K(Y,G)/(\lambda_{-1}(\N_Y X))
    \]
is induced by $i^* \colon \K(X,G) \to \K(Y,G)$.
\end{enumeratei}
\end{proposition}

\begin{proof}
 From the self-intersection formula (Theorem~\ref{thm:selfintersection}) we
see that the composition $i^*i_* \colon \K(Y,G) \to \K(Y,G)$ is
multiplication by $\lambda_{-1}(\N_Y^\vee X)$, so $i_*$ is injective. We
get part~\refpart{prop:K-rigid->}{1} from this and from the localization
sequence.

Part~\refpart{prop:K-rigid->}{2} follows easily from
part~\refpart{prop:K-rigid->}{1}, together with the following elementary
fact.

\begin{lemma}
Let $A$, $B$ and $C$ be rings, $f \colon B \to A$ and $g \colon B \to C$
ring homomorphisms. Suppose that there exist a homomorphism of abelian
groups $\phi\colon A \to B$ such that:
\begin{enumeratei}
\item The sequence
    \[
    0 \longrightarrow A \overset{\phi}\longrightarrow B \overset g
    \longrightarrow C \longrightarrow 0.
    \]
is exact;

\item the composition $f\circ \phi \colon A \to A$ is the multiplication
by a central element $a \in A$ which is not a zero divisor.
\end{enumeratei}
Then $f$ and $g$ induce an isomorphism of rings
    \[
    (f,g) \colon B \to A
    \displaytimes_{A/(a)} C,
    \]
where the homomorphism $A \to A/(a)$ is the projection,
and the one $C \to A/(a)$ is induced by the isomorphism $C \simeq B/\im \phi$
  and the projection $g \colon B \to A$.\qedhere
\end{lemma}
\end{proof}

\subsection{The theorem of reconstruction from the strata}
\label{sec:maintheorem}

This section is entirely dedicated to the proof of our main theorem. Let
us recall what it says. Let $G$ act on $X$ with the usual hypotheses.
Consider the strata $X_s$ defined in
Subsection~\ref{subsec:stratdimension}, and the specialization
homomorphisms
    \[
    \sp_{X,s}^t \colon \K(X_{t},G) \longrightarrow \K\bigl(\N_{s,t},
    G\bigr)
    \]
defined in
Subsection~\ref{subsec:specializationsnormal}.

\begin{theorem}[The theorem of reconstruction from the
strata]\call{thm:maintheorem}
The homomorphism
    \[
    \K(X,G) \longrightarrow \prod_{s=0}^n \K(X_s,G)
    \]
obtained from the restrictions $\K(X,G) \to \K(X_s,G)$ is injective. Its
image consists of the sequences $(\alpha_s) \in \prod_{s=0}^n \K(X_s,G)$
with the property that for each $s= 1$, \dots,~$n$ the pullback of $\alpha
\in  \K(X_s,G)$ to $\K\bigl(\N_{s,s-1}, G\bigr)$ coincides with\/
$\sp_{X,s}^{s-1}(\alpha_{s-1}) \in \K\bigl(\N_{s,s-1}, G\bigr)$.
\end{theorem}

In other words, we can view $\K(X,G)$ as a fiber product
    \begin{align*}
    \K(X,G) \simeq {}&\K(X_n,G) \displaytimes_{\K(\N_{n,n-1},G)}
    \K(X_{n-1},G)
    \displaytimes_{\K(\N_{n-1,n-2},G)}\\
    &\quad\ldots
    \displaytimes_{\K(\N_{2,1},G)}  \K(X_1,G) \displaytimes_{\K(\N_{1,0},G)}
    \K(X_0,G).
    \end{align*}

Here is our starting point.

\begin{proposition}
$X_s$ is K-rigid in $X$.
\end{proposition}

\begin{proof} This follows from Proposition
\refall{prop:firststratification}{2}, and
Lemma~\ref{lem:criterionK-rigid}.
\end{proof}

So from  Proposition \refall{prop:K-rigid->}{2} applied to the closed
embedding $i_s \colon X_s \into X_{\le s}$, we get an isomorphism
   \[
   \K(X_{\le s},G) \simeq \K(X_s,G)
   \displaytimes_{\K(X_s,G)/(\lambda_{-1}(\N_s^\vee))}
   \K(X_{<s},G).
   \]

We can improve on this.

\begin{proposition}\call{prop:fiberproduct-1stratum} The restrictions
    \[
    \K(X_{\le s},G) \longrightarrow \K(X_s,G)\quad \text{and} \quad
    \K(X_{\le s},G)\longrightarrow \K(X_{<s},G)
    \]
induce an isomorphism
    \[
    \K(X_{\le s},G) \larrowsim \K(X_s,G)
    \displaytimes_{\K(\N_s^0,G)} \K(X_{<s},G),
    \]
where the homomorphism $\K(X_s,G) \to \K(\N_s^0,G)$ is the pullback, while
    \[
    \sp_{Y,s} \colon \K(X_{<s},G) \to \K(\N_s^0,G)
    \]
is the specialization.
\end{proposition}

\begin{proof} Let us start with a lemma.

\begin{lemma}\label{lem:zerosection}
The restriction homomorphism
$\K(X_s,G) \to \K(\N_s^0,G)$ is surjective, and its kernel is the
ideal $\bigl(\lambda_{-1}(\N_s^\vee)\bigr) \subseteq \K(X_s,G)$.
\end{lemma}
\begin{proof} Since the complement of the zero
section
$s_0 \colon X_s
\into \N_s$ coincides with $(\N_s)_{<s}$
(Proposition \refall{prop:firststratification}{2}), we can apply
Proposition \refall{prop:K-rigid->}{1}
to the normal bundle $\N_s$, and conclude that there is an exact sequence
    \[
    0 \longrightarrow \K(X_s,G) \overset{s_{0*}}\longrightarrow \K(\N_s,G)
    \longrightarrow \K(\N_s^0,G) \longrightarrow 0.
    \]

Now, $s_0^* \colon \K(\N_s,G) \to \K(X_s,G)$ is an isomorphism, and the
composition $s_0^*s_{0*} \colon \K(X_s,G) \to \K(X_s,G)$ is multiplication
by $\lambda_{-1}(\N_s^\vee)$, because of the self-intersection
formula~\ref{thm:selfintersection}, and this implies the thesis.
\end{proof}

Therefore the restriction homomorphism $\K(X_s,G) \to
\K(\N_s^0,G)$ induces an isomorphism of
$\K(X_s,G)/\bigl(\lambda_{-1}(\N_s^\vee)\bigr)$ with $\K(\N_s^0,G)$; the
Proposition follows from this, and from
Proposition~\ref{prop:compspecializations}.
\end{proof}

Now we proceed by induction on the largest integer $s$ such that $X_s
\neq \emptyset$. If $s = 0$ there is nothing to prove. If $s > 0$, by
induction hypothesis the homomorphism
    \[
    \K(X_{<s},G) \longrightarrow \K(X_{s-1},G)
\displaytimes_{\K(\N_{s,s-1},G)} \dots \displaytimes_{\K(\N_{1,0},G)}
\K(X_0,G)
    \]
induced by restrictions is an isomorphism; so from
Proposition~\ref{prop:fiberproduct-1stratum} we see that to
prove Theorem~\ref{thm:maintheorem} it is sufficient to show that if
$\alpha_s \in \K(X_s,G)$, $\alpha_{<s} \in \K(X<s,G)$, $\alpha_{s-1}$
is the restriction of $\alpha_{<s}$ to $\K(X_{s-1},G)$, $\alpha_s^0$ is
the pullback of $\alpha_s$ to $\K(\N_s^0,G)$ and $\alpha_{s,s-1}$ is
the pullback of $\alpha_s$ to $\K(\N_{s,s-1},G)$, then
$\sp_{X,s}(\alpha_{<s}) = \alpha_s^0$ if and only if
$\sp_{X,s}^{s-1}(\alpha_{s-1}) = \alpha_{s,s-1}$. But the diagram
    \[
    \xymatrix{
    {}\K(X_{<s},G)\ar[r]^{\sp_{X,s}}\ar[d]& {}\K(\N_s^0,G)\ar[d]\\
    {}\K(X_{s-1},G)\ar[r]^{\sp_{X,s}^{s-1}}&{}\K(\N_{s,s-1},G)
    },
    \]
where the colums are restriction homomorphisms, is commutative
(Proposition~\ref{prop:compspecializations}); hence it suffices
to show that the restriction homomorphism
    \[
    \K(\N_s^0,G) \longrightarrow \K(\N_{s,s-1},G)
    \]
is injective. To prove this we may suppose that the
action of $G$ on $X_s$ is connected, that is, $X_s$ is not a nontrivial
disjoint union of open invariant subspaces. In this case the toral component
of the isotropy group of a point of $X_s$ is constant.

Set $E = \N_s$, and consider the eigenspace decomposition $E =
\bigoplus_{\chi
\in \widehat T} E_\chi$. We obtain a decomposition $E = \bigoplus_i E_i$
by grouping together $E_\chi$ and $E_{\chi'}$ when the characters $\chi$
and $\chi'$ are multiple of a common primitive character in $\widehat T$.
Then clearly a geometric point of $E$ is in $E_{s-1}$ if and only if
exactly one of its components according to the decomposition $E =
\bigoplus_i E_i$ above is nonzero. In other words, $\N_{s,s-1}$ is the
disjoint union
$\coprod_i E_i^0$, where $E_i^0$ is embedded in $E$ by setting all the
other components equal to 0. The same argument as in the proof of
Lemma~\ref{lem:zerosection} shows that the kernel of the pullback
$\K(X_s, G) \to \K(E_i^0,G)$ is generated by $\lambda_{-1} E_i$, so the
kernel of the pullback
    \[
    \K(X_s,G) \longrightarrow \K(\N_{s,s-1},G) = \bigoplus_i\K(E_i^0,G)
    \]
equals $\cap_i
(\lambda_{-1}E_i)$; hence we need to show that
$(\lambda_{-1} E) = \cap_i (\lambda_{-1} E_i)$.

This is done as follows. Choose a splitting $G = D\times T$: we have
$\K(X_s,G) = \K(X_s,D) \otimes \R T$, as in the beginning of the proof of
Lemma~\ref{lem:criterionK-rigid}. First of all, we have
$\lambda_{-1}E = \prod_i
\lambda_{-1} E_i$. Furthermore, for each $i$ we can choose a primitive
character $\chi_i$ in $\widehat T$ such that all the characters which
appear in the decomposition of $E_i$ are multiples of $\chi_i$; from this we
see that $\lambda_{-1}E_i$ is of the form $\sum_{k=m_i}^{n_i} r_{i,k}
\chi_i^k$, where $m_i$ and $n_i$ are (possibly negative) integers, $r_{i,k}
\in
\K[0](X_s,G)$, $r_{i,m}$ and $r_{i,n}$ are invertible. Then the conclusion
of the proof follows from the following fact.

\begin{lemma}\label{lem:int=prod}
Let $A$ be a ring, $H$ a free finitely generated abelian group, $\chi_1$,
\dots,~$\chi_r$ linearly independent elements of $H$. Let
$\gamma_1$, \dots,~$\gamma_r$ be elements of the group ring $AH$ of the
form $\gamma_i = \sum_{k=m_i}^{n_i} r_{i,k} \chi_i^k$, where the $r_{i,k}$
are central elements of $A$ such that $r_{i,m_i}$ and $r_{i,n_i}$ are
invertible. Then we have an equality of ideals $(\gamma_1 \ldots \gamma_r)
= (\gamma_1)\cap \ldots\cap(\gamma_r)$ in $AH$.
\end{lemma}

\begin{proof}
By multiplying each $\gamma _i$ by
$r_{i,m_i}^{-1}\chi _i^{-m_i}$ we may assume that $ \gamma _i$ has the form
$1+a_{i,1}\chi _i+\cdots +a_{i,s_i}\chi _i^{s_i}$ with $s_i\geq 0$ and
$a_{i,s_i}$ a central unit in $A$. We will show that for any
$i\neq j$ the relation
$\gamma _i\mid q\gamma _j$, $q\in AH$ implies $\gamma _i\mid q$; from
this the thesis follows with a straightforward induction. We may assume
that $r =2$, $i = 1$ and $j = 2$.

Since $\chi _1$, $\chi_2$ are $\mathbb{Z}$-linearly independent elements
of $H$, we may complete them to a maximal
$\mathbb{Z}$-independent sequence $\chi _1,\dots ,~\chi _n$ of $H$; this
sequence generates a subgroup $H' \subseteq H$ of finite index.

Suppose at first that $H' = H$, so that
$AH = A\bigl[\chi _1^{\pm 1},...,\chi _n^{\pm 1}\bigr]$. Replacing $A$ by
$A\bigl[\chi _3^{\pm 1}, \dots, \chi_n^{\pm 1}\bigr]$, we may assume that
$AH = A\bigl[\chi _1^{\pm 1},\chi _2^{\pm 1}\bigr]$.

If $p\gamma _1=q\gamma _2$, we can multiply this equality by a
sufficiently high power of $\chi_1 \chi_2$ and assume that $p$ and $q$ are
polynomials in $A[\chi_1, \chi_2]$. Since $\gamma_2$ is a polynomial in
$A[\chi_2]$ with central coefficients and invertible leading coefficient,
the usual division algorithm allows us to write $p = s \gamma_2 + r \in
A[\chi_1, \chi_2]$, where $r$ is a polynomial whose degree in $\chi_2$ is
less than $s_2 = \deg_{\chi_2}\gamma_2$. By comparing the degrees in
$\chi_2$ in the equality $r\gamma_1 = (q - s\gamma_1)\gamma_2$ we see
that $q - s\gamma_1$ must be zero, and this proves the result.

In the general case, choose representatives $u_1$, \dots,~$u_r$ for the
cosets of $H'$ in $H$; then any element $f$ of $AH$ can be written
uniquely as $\sum_{i=1}^r u_i f_i$ with $f_i \in AH'$. Then from the
equality $\bigl(\sum_i u_i p_i\bigr) \gamma_1 = \sum_i u_i q_i$ we get
$p_i \gamma_1 = q_i \gamma_2$ for all $i$, because $\gamma_1$ and
$\gamma_2$ are in $AH'$; hence the thesis follows from the previous case.
\end{proof}

\section{Actions with enough limits}\label{sec:limits}

Let us start with some preliminaries in commutative algebra.

\subsection{Sufficiently deep modules}

Let $A$ be a finitely generated flat \cm $\mathbb{Z}$-algebra, such that each
of the fibers of the morphism $\spec A \to \spec \mathbb{Z}$ has pure
dimension
$n$. If $V$ is a closed subset of $\spec A$, we define the \emph{fiber
dimension} of $V$ to be the largest of the dimensions of the fibers of $V$
over $\spec \mathbb{Z}$, and its
\emph{fiber codimension} to be $n$ minus its fiber dimension. We say that
$V$ has \emph{pure fiber dimension} if all the fibers of $V$ have the same
fiber dimension at all points of $V$ (of course some of the fibers may be
empty).

The fiber dimension and codimension of an ideal in $A$ will be the fiber
dimension and codimension of the corresponding closed subset of $\spec A$.

\begin{definition}\call{def:suff-deep}
Let $M$ be an $A$-module. Then we say that $M$ is
\emph{sufficiently deep} if the following two conditions are satisfied.

\begin{enumeratei}

\itemref{1}  All associated primes of M have fiber codimension 0.

\item $\ext_A^1(A/\mathfrak{p}, M) = 0$ for all primes $\mathfrak{p}$ in
$A$ of fiber codimension at least 2.

\end{enumeratei}
\end{definition}

Here are the properties that we need.

\begin{proposition}\call{prop:suff-deep}\hfil
\begin{enumeratei}

\itemref{1}  If\/ $0 \to M' \to M \to M'' \to
0$ is an exact sequence of $A$-modules, $M'$ and $M''$ are sufficiently
deep, then
$M$ is sufficiently deep.

\itemref{2} Direct limits and direct sums of
sufficiently deep modules are sufficiently deep.

\itemref{3}  If $N$ is an abelian group, then
$N\otimes_\mathbb{Z} A$ is sufficiently deep.

\itemref{4} If $M$ is a sufficiently deep
$A$-module, then\/
$\ext_A^1(N,M) = 0$ for all
$A$-modules $N$ whose support has fiber codimension at least 2.
\end{enumeratei}
\end{proposition}

\begin{proof} Part~\refpart{prop:suff-deep}{1} is obvious.

Part~\refpart{prop:suff-deep}{2} follows from the fact that $A$ is
noetherian, so formation of $\ext_A^1(A/\mathfrak{p}, -)$ commutes with
direct sums and direct limits.

Let us prove part~\refpart{prop:suff-deep}{3}. From
part~\refpart{prop:suff-deep}{2} we see that we may assume that $N$ is
cyclic. If $N = \mathbb{Z}$, then $M = A$,  and the statement follows from
the facts that $A$ is \cm, and that the height of a prime
ideal is at least equal to its fiber codimension.

Assume that $N = \mathbb{Z}/m\mathbb{Z}$, so that $M = A/mA$. The
associated primes of $M$ are the generic points of the fibers of $A$ over
the primes dividing $m$, so condition \refall{def:suff-deep}{1} is
satisfied.

Take a prime $\mathfrak{p}$ of $A$ of
fiber codimension at least 2, and consider the exact sequence
    \[
    0 = \ext_A^1(A/\mathfrak{p}, A) \to \ext_A^1(A/\mathfrak{p},
    A/mA)
    \to \ext_A^2(A/\mathfrak{p}, A)
    \overset m{\to} \ext_A^2(A/\mathfrak{p}, A).
    \]
If the characteristic of $A/\mathfrak{p}$ is positive, then the height of
$\mathfrak{p}$ is at least~3, so $\ext_A^2(A/\mathfrak{p}, A) =0$,
because $A$ is \cm, and we are done. Otherwise, we have an
exact sequence
    \[
    0 \longrightarrow A/\mathfrak{p}
    \overset m{\longrightarrow} A/\mathfrak{p} \longrightarrow A/\bigl((m) +
    \mathfrak{p}\bigr)\longrightarrow 0;
    \]
but the height of $(m) + \mathfrak{p}$ is at least 3, so
$\ext_A^2\bigl(A/\bigl((m) + \mathfrak{p} \bigr), A \bigr) = 0$. From
this we deduce that multiplication by $m$ is injective on
$\ext_A^2(A/\mathfrak{p}, A)$, and this concludes the proof of
part~\refpart{prop:suff-deep}{3}.

For part~\refpart{prop:suff-deep}{4}, notice first of all that if $N$
is a finitely generated
$A$-module of fiber codimension at least 2 then we can filter $N$ with
successive quotients of type $A/\mathfrak{p}$, where $\mathfrak{p}$ is a
prime of fiber codimension at least 2, so $\ext_A^1(N,M) = 0$.

If $N$ is not finitely generated and $0 \to M \to E \to N \to 0$ is an
exact sequence of $A$-modules, $N'$ is a finitely generated submodule of
$N$, and $E'$ is the pullback of $E$ to $N'$, then the sequence $0 \to M
\to E' \to N'$ splits; but because of part~\refpart{prop:suff-deep}{1} of
the definition we have $\hom_A(N',M) = 0$, hence there is a unique copy
of $N'$ inside
$E'$. Hence there is a unique copy of $N$ inside $E$, and the sequence
splits. This completes the proof of the Proposition. \end{proof}

\subsection{Sufficiently deep actions}

Let $G$ be a diagonalizable group scheme of finite type
over $S$; all the actions will
be upon noetherian separated regular algebraic spaces, as in our setup.
The ring of representations $\R G = \mathbb{Z} \widehat G$ is a finitely
generated flat \cm
$\mathbb{Z}$-algebra, and each of the fibers of the morphism $\spec \R G \to
\spec \mathbb{Z} $ has pure dimension equal to the dimension of $G$.

\begin{definition} We say that the action of $G$ on $X$ is
\emph{sufficiently deep} when the $\mathrm{R}G$-module $\K(X,G)$
is sufficiently deep.
\end{definition}

\begin{theorem}\call{thm:suff-deep->} Suppose that a diagonalizable group
scheme of finite type
$G$ acts on a noetherian regular separated algebraic space $X$, and that
the action is sufficiently deep. Then the restriction homomorphism
$\K(X,G) \to \K(X^{G_0},G)$ is injective,
and its image is the intersection of the images of the restriction
homomorphisms\/ $\K(X^T,G)\to \K(X^{G_0},G)$, where
$T$ ranges over all subtori of $G$ of codimension 1.
\end{theorem}

\begin{proof} We need some preliminaries.

\begin{lemma}\call{lem} Suppose that $G$ acts on $X$ with
stabilizers of constant dimension $s$. Then the support of
$\K(X,G)$ as an
$\mathrm{R}G$-module has pure fiber dimension $s$, and any associated
prime of $\K(X,G)$ has pure fiber dimension $s$.
\end{lemma}

\begin{proof} Suppose first of all that $s$ is 0. Then  it follows easily
from Thomason's localization theorem (\cite{th3}) that the support of
$\K(X,G)$ has fiber dimension 0, and from this we see that every
associated prime must have fiber dimension 0.

In the general case, we may assume that the action is connected (that is,
$X$ is not a nontrivial disjoint union of open invariant subspaces); then
there will be a splitting $G = H \times_S T$, where $H$ is a diagonalizable
group scheme of finite type acting on
$X$ with finite stabilizers, and $T$ is a totally split torus that acts
trivially on
$X$. In this case
    \[
    \K(X,G) = \K(X,H)\otimes_\mathbb{Z} RT
    = \K(X,H)\otimes_{\mathrm{R}H} \mathrm{R}G.
    \]
The proof is concluded by applying the following lemma.

\begin{lemma} Let $A$ be a flat \cm
$\mathbb{Z}$-algebra of finite type, $A\to B$ a smooth homomorphism of
finite type with fibers of pure dimension $s$. Suppose that $M$ is an
$A$-module whose support has fiber dimension 0. Then $M\otimes_A B$ has
support of pure dimension $s$, and each of its associated primes has fiber
dimension
$s$.
\end{lemma}

\begin{proof} Since tensor product commutes with taking direct limits and
$B$ is flat over $A$, we may assume that $M$ is of finite type over $A$.
By an obvious filtration argument, we may assume that $M$ is of the form
$A/\mathfrak{p}$, where $\mathfrak{p}$ is a prime ideal of fiber dimension
0. In this case the only associated primes of $M \otimes_A B$ are
the generic components of the fiber of $\spec B$ over $\mathfrak{p}$, and
this proves the result. \end{proof}

\noqed \end{proof}

\begin{lemma}\call{lem:unmixing} Suppose that $X$ and $Y$ are
algebraic spaces on which $G$ acts with stabilizers of constant dimension
respectively $s$ and $t$. If $N$ is an $\mathrm{R}G$-submodule of
$\K(Y,G)$ and $t<s$, then there is no nontrivial homomorphism of $\R
G$-modules from $N$ to
$\K(X,G)$.
\end{lemma}

\begin{proof} Given such a nontrivial homomorphism $N \to \K(X,G)$, call
$I$ its image. The support of $I$ has fiber dimension at most $t$, so
there is an associated prime of fiber dimension at most $t$ in
$\K(X,G)$, contradicting Lemma~\ref{lem}.
\end{proof}

Now we prove Theorem~\ref{thm:suff-deep->}. Let $n$ be the
dimension of
$G$, so that $X_n = X^{G_0}$. First of all, let us show that the natural
projection
    \[
    \K(X,G) \longrightarrow \K(X_n,G) \displaytimes_{\K(\N_{n,n-1}, G)}
    \K(X_{n-1},G)
    \]
is an isomorphism. This will be achieved by showing
that for all $s$ with $0 \le s \le n-1$ the natural projection
$\K(X,G) \longrightarrow P_s$ is an isomorphism, where we have set
    \begin{align*}
    P_s =  {}&\K(X_n,G) \displaytimes_{\K(\N_{n,n-1}, G)}
    \K(X_{n-1},G) \displaytimes_{\K(\N_{n-1,n-2},G)}\\
    &\quad\cdots
    \displaytimes_{\K(\N_{s+1,s},G)} \K(X_s,G).
    \end{align*}

For $s = 0$ this is our main Theorem~\ref{thm:maintheorem}, so we proceed
by induction. If
$s < n-1$ and the projection above is an isomorphism, we have an exact
sequence
    \[
    0 \longrightarrow \K(X,G) \longrightarrow P_{s+1} \times
    \K(X_s,G) \longrightarrow \K(\N_{s+1,s},G),
    \]
where the last arrow is the difference of the composition of the
projection $P_{s+1} \to \K(X_{s+1}, G)$ with the pullback $\K(X_{s+1}, G)
\to \K(\N_{s+1,s},G)$, and of the specialization homomorphism $\K(X_s,G) \to
\K(\N_{s+1,s},G)$. If we call $N$ the image of this difference, we
have an exact sequence of $\mathrm{R}G$-modules
    \[
    0 \longrightarrow \K(X,G) \longrightarrow P_{s+1} \times \K(X_s,G)
    \longrightarrow N \to 0,
    \]
and the support of $N$ is of fiber dimension at most $s \le n-2$ by
Lemma~\ref{lem}, hence it is of fiber codimension at least 2.
It follows from the fact that $\K(X,G)$ is sufficiently deep and from
Proposition \refall{prop:suff-deep}{4} that this sequence splits. From
the fact that $\K(X,G)$ has only associated primes of fiber dimension~0 we
see that the pullback map
$\K(X_s,G) \to N$ must be injective, and from Lemma~\ref{lem:unmixing} that a
copy of
$N$ living inside $ P_{s+1} \times \K(X_s,G)$ must in fact be contained
in $\K(X_s,G)$; this implies that the projection $\K(X,G)
\to P_{s+1}$ is an isomorphism.

So the projection
    \[
    \K(X,G) \longrightarrow \K(X_n,G) \times_{\K(\N_{n,n-1}, G)}
    \K(X_{n-1},G)
    \]
is an isomorphism. Then the kernel of the specialization
homomorphism from $\K(X_{n-1},G)$ to $\K(\N_{n,n-1}, G)$ maps
injectively in $\K(X,G)$, so it must be 0, again because $\K(X,G)$ has only
associated primes of fiber dimension~0. Furthermore
$X_{n-1}$ is the disjoint union of the $X_{n-1}^T$ when $T$ ranges over all
finite subtori of $G$ of codimension 1, and similarly for $\N_{n,n-1}$.
On the other hand, because of our main theorem applied to the action
of $G$ on $X^T$, we have the natural isomorphism $\K(X^T,G) \to
\K(X^{G_0}, G) \times_{\K(\N_{n,n-1}^T,G)}\K(X_{n-1}^T,G)$, and this
completes the proof of Theorem~\ref{thm:suff-deep->}. \end{proof}

\subsection{Actions with enough limits are sufficiently deep}

For the rest of this section $S$ will be the spectrum of a field $k$, $G$
is a diagonalizable group scheme  of finite type acting on a smooth
separated scheme
$X$ of finite type over $k$; call $M$ the group of one-parameter subgroups
$\gm[k]
\to G$ of
$G$. There is a natural Zariski topology on $M \simeq \mathbb{Z}^n$ in which
the closed subsets are the loci of zeros of  sets of polynomials in the
symmetric algebra $\mathop{\mathrm{Sym}}_{\mathbb{Z}}^{\bullet} M^\vee$;
we refer to this as the \emph{Zariski topology on $M$}.

We will denote, as usual, by $G_0$ the toral component of $G$. If $n$ is
the dimension of $G$, then $X_n$ = $X^{G_0}$. Furthermore, if we choose a
splitting $G \simeq G_0 \times G/G_0$ we obtain an isomorphism of rings
    \[
    \K(X^{G_0},G) \simeq \K(X^{G_0}, G/G_0) \otimes \R G_0
    \]
(\cite[Lemme~5.6]{th4}).

\begin{definition}\call{def:admitslimits} Suppose that $k$ is
algebraically closed. Consider a one parameter subgroup $H = \gm[k] \to G$,
with the corresponding action of $\gm[k]$ on $X$. We say that this one
parameter subgroup \emph{admits limits} if for every closed point $x \in
X$, the morphism $\gm[k] \to X$ which sends $t \in G$ to $tx$ extends to a
morphism
$\mathbb{A}^1 \to X$. The image of $0 \in \mathbb{A}^1(k)$ in $X$ is
called \emph{the limit of $x$ for the one parameter subgroup $H$}.

We say that the action of $G$ on $X$ \emph{admits enough limits}
if the one parameter subgroups of $G$ which admit limits form a
Zariski-dense subset of $M$.

If $k$ is not algebraically closed, then we say that the action admits
enough limits if the action obtained after base change to the algebraic
closure of $k$ does. \end{definition}

\begin{remark}
One can show that the locus of $1$-parameter subgroups of $G$ admitting
limits is defined by linear inequalities, so the definition can be stated
in more down to earth terms (we are grateful to the referee for pointing
this out).
\end{remark}

The notion of action with enough limits is a weakening of the notion of
\emph{filtrable} action due to Brion. More precisely, an action has enough
limits if it satisfies condition~(i) in \cite[Definition~3.2]{brion1}; there
is also a condition~(ii) on closures of strata.

The main case when the action admits enough limits is when $X$
is complete; in this case of course every one-parameter subgroup admits
limits. Another case is when the action of $G_0 = \gm[k]^n$ on $X$ extends
to an action of the multiplicative monoid $\mathbb{A}^n$.
Also, we give a
characterization of toric varieties with enough limits in
Prop. ~\ref{charact}.

\begin{theorem}\call{thm:->suff-deep} Suppose that a diagonalizable group
scheme of finite type
$G$ over a perfect field $k$ acts on a smooth separated scheme of finite
type $X$ over $k$. If the action of $G$ admits enough limits, then it
is sufficiently deep.
\end{theorem}

By putting this together with Theorem~\ref{thm:suff-deep->} we get the
following.

\begin{corollary}\call{cor:maincorollary}
Suppose that a diagonalizable group scheme  of finite type
$G$ over a perfect field $k$ acts on a smooth separated scheme of finite
type $X$ over $k$. If the action of $G$ admits enough limits, then the
restriction homomorphism
    \[
    \K(X,G) \longrightarrow \K(X^{G_0},G)
    \]
is injective, and its image is the intersection of
the images of the restriction homomorphisms\/ $\K(X^T,G)\to
\K(X^{G_0},G)$, where $T$ ranges over all subtori of $G$ of codimension 1.
\end{corollary}

For example, consider the following situation, completely analogous to the
one considered in \cite[Corollary~7]{brion2} and in \cite{gkmp}. Let $G$ be
an
$n$-dimensional torus acting on a smooth complete variety
$X$ over an algebraically closed field $k$. Assume that the fixed point set
$X^{G_0} = X_n$ is zero-dimensional, while $X_{n-1}$ is 1-dimensional.
Set $X^{G_0} = \{x_1, \dots, x_t\}$, and call $P_1$, \dots,~$P_r$ the
closures in $X$ of the connected components of $X_{n-1}$. Then each $P_j$
is isomorphic to $\mathbb{P}^1$, and contains precisely two of the fixed
points, say $x_i$ and $x_{i'}$. Call $D_j$ the kernel of the action of
$G$ on $P_j$; then the image of the restriction
homomorphism
    \[
    \K(P_j,G) \to \K(x_i,G) \times \K(x_{i'},G) = \K(k)\otimes \R
    G \times \K(k)\otimes \R G
    \]
consists of the pairs of elements
    \[
    ( \alpha,
    \beta) \in \K(k)\otimes \R G \times \K(k)\otimes \R G
    \]
whose images in
$\K(k) \otimes \R D_j$ coincide (this follows immediately from
Theorem~\ref{thm:maintheorem}). From this and from
Corollary~\ref{cor:maincorollary} we get the following.

\begin{corollary}\label{cor:generic}
In the situation above, the restriction map
    \[
    \K(X,G) \longrightarrow \prod_{i=1}^t
    \K(x_i,G) = \prod_{i=1}^n \K(k) \otimes \R G
    \]
is injective. Its image
consist of all elements $(\alpha_i)$ such that if $x_i$ and $x_{i'}$ are
contained in some $P_j$, then the restrictions of $\alpha_i$ and
$\alpha_{i'}$ to $\K(k) \otimes \R D_j$ coincide.
\end{corollary}

Theorem~\ref{thm:->suff-deep} is proved in the next subsection.

\subsection{Bia\l ynicki-Birula stratifications}

Let us prove Theorem~\ref{thm:->suff-deep}: like in \cite{brion1}, the
idea is to use a Bia\l ynicki-Birula stratification. We will prove the
following.

\begin{proposition}\label{prop:describesufflimits}
Suppose that a diagonalizable group scheme  of finite type
$G$ over a perfect field $k$ acts with enough limits on a smooth separated
scheme of finite type $X$ over $k$. Then the $\R G$-module
$\K(X,G)$ is obtained by taking finitely many successive
extensions of $\R G$-modules of the form $N \otimes_\mathbb{Z}\R G$,
where $N$ is an abelian group.
\end{proposition}

Theorem~\ref{thm:->suff-deep} follows from this, in view of
Proposition~\ref{prop:suff-deep}, parts \refpart{prop:suff-deep}{1} and
\refpart{prop:suff-deep}{3}. Let us prove the Proposition.

First of all, let us
assume that $k$ is algebraically closed. We will only consider closed
points, and write $X$ for $X(k)$.

It is a standard fact that the one-parameter subgroups
$H = \gm[k] \to G_0$ with the property that $X^{G_0} = X^H$ form a nonempty
Zariski open subset of $M$, so we can choose one with this property that
admits enough limits. There is a (discontinuous) function $X \to X^{G_0}$
sending each point to its limit. Let $T_1$, \dots,~$T_s$ be the connected
components of $X^{G_0}$; call
$E_i$ the inverse image of $T_i$ in $X$, and $\pi_j\colon E_j
\to T_j$ the restriction of the limit function. The following is a
fundamental result of Bia\l ynicki-Birula.

\begin{theorem}[Bia\l ynicki-Birula] \call{thm:BB} In the situation above:

\begin{enumeratei}

\item The
$E_j$ are smooth locally closed $G$-invariant subvarieties of $X$.

\item  The functions $\pi_j\colon E_j \to T_j$ are
$G$-invariant morphisms.

\itemref{1}  For each $j$ there is a representation
$V_j$ of $H$ and an open cover $\{U_\alpha\}$ of $T_j$, together with
equivariant isomorphisms $\pi_j^{-1}(U_\alpha) \simeq U_\alpha \times V_j$,
such that the restriction $\pi_j\colon  \pi_j^{-1}(U_\alpha) \to U_\alpha$
corresponds to the projection $U_\alpha \times V_j \to U_\alpha$.

\item If $x$ is a point of $T_j$, then the normal space to $E_j$
in $X$ at $x$ is the sum of the negative eigenspaces in the tangent space
to $X$ at $x$ under the action of $H$.

\end{enumeratei}

\end{theorem}

Of course in part~\refpart{thm:BB}{1} we may take $V_j$ to be
the normal bundle to $E_j$ in $X$ at any point of $T_j$.

This theorem is proved in
\cite{b-b}; we should notice that the condition that $X$ is covered by
open invariant quasiaffine subsets is always verified, thanks to a result
of Sumihiro (\cite{sumihiro}).

Now we remove the hypothesis that $k$ is algebraically closed: here is the
variant of Bia\l ynicki-Birula's theorem that we need.

\begin{theorem}\call{thm:BBvariant} Suppose that a diagonalizable group
scheme of finite type
$G$ over a perfect field $k$ acts with enough limits on a smooth separated
scheme of finite type $X$ over $k$. Let $Y_1$, \dots,~$Y_r$ be the
connected components of $X^{G_0}$; there exists a stratification $X_1$,
\dots,~$X_r$ of
$X$ in locally closed $G$-invariant smooth subvarieties, together with
$G$-equivariant morphisms $\rho_i\colon X_i \to Y_i$, such that:
\begin{enumeratei}

\itemref{1} $X_i$ contains $Y_i$ for all $i$, and
the restriction of
$\rho_i$ to $Y_i$ is the identity.

\itemref{2} If $U$ is an open affine subset of
$Y_i$ and\/
$\N_U$ is the restriction of the normal bundle\/
$\N_{Y_i}{X_i}$ to $U$, then there is a $G$-equivariant isomorphism
$\rho_i^{-1}(U) \simeq \N_U$ of schemes over $U$.

\itemref{3} In the eigenspace
decomposition of the restriction of
$\N_{X_i}{X}$ to $Y_i$, the subbundle corresponding to the
trivial character of $G_0$ is 0.
\end{enumeratei}
\end{theorem}

\begin{proof} Let $\overline X = X \times_{\spec k} \spec \overline 
k$, and call
$\Gamma$ the Galois group of $\overline k$ over $k$. Choose a one
parameter subgroup $H = \gm[k] \to G$ as before. Let $T_1$,
\dots,~$T_s$ be the connected components of $\overline X^{G_0}$,
$\pi_j\colon E_j \to T_j$ as in Bia\l ynicki-Birula's theorem. The $Y_i$
correspond to the orbits of the action of $\Gamma$ on $\{T_1, \ldots,
T_s\}$; obviously
$\Gamma$ also permutes the $E_j$, so we let $X_1$, \dots,~$X_r$ be the
smooth subvarieties of $X$ corresponding to the orbits of $\Gamma$ on
$\{E_1, \ldots, E_s\}$. The $\pi_j\colon E_j \to T_j$ descend to morphisms
$X_i \to Y_i$. Properties \refpart{thm:BBvariant}{1}
and \refpart{thm:BBvariant}{3} are obviously satisfied,
because they are satisfied after passing to
$\overline k$.

We have to prove \refpart{thm:BBvariant}{2}. Let $E$ be the inverse
image of $U$ in $X_i$,
$I$ the ideal of $U$ in the algebra $k[E]$. Because $U$ is affine, $I/I^2$
is a projective $k[U]$-module, and $G$ is diagonalizable, the projection $I
\to I/I^2$ has a $k[U]$-linear and $G$-equivariant section
$I/I^2 \to I$. This induces a
$G$-equivariant morphism of $U$-schemes $E \to \N_U$, sending $U$ to
the 0-section, whose differential at the zero section is the identity
(notice that $\N_U$ is also the restriction to $U$ of the relative
tangent bundle ${\rm T}_{X/U}$). We want to show that this is an
isomorphism; it is enough to check that this is true on the fibers,
so, let $V$ be one of the fibers of $\N_U$ on some point $p \in U$.
According to part~\refpart{thm:BBvariant}{2} of the theorem of Bia\l
ynicki-Birula, the fiber of
$X$ on $p$ is
$H$-equivariantly isomorphic to $V$; hence an application of the
following elementary lemma concludes the proof of
Theorem~\ref{thm:BBvariant}.

\begin{lemma} Suppose that $\gm[k]$ acts linearly with positive weights on
a finite dimensional vector space $V$ over a field $k$. If $f \colon V \to
V$ is an equivariant polynomial map whose differential at the origin is an
isomorphism, then $f$ is also an isomorphism.
\end{lemma}
\begin{proof}
First all, notice that because of the positivity of the weights, we have
$f(0) = 0$. By composing
$f$ with the inverse of the differential of
$f$ at the origin, we may assume that this differential of $f$ is the
identity. Consider the eigenspace decomposition $V = V_1 \oplus V_2 \oplus
\dots \oplus V_r$, where $\gm[k]$ acts on $V_i$ with a character $t \mapsto
t^{m_i}$, and $0 < m_1 < m_2 < \dots < m_r$. Choose a basis of
eigenvectors of $V$; we will use groups of coordinates $x_1$,
\dots,~$x_r$, where $x_i$ represents the group of elements of the dual
basis corresponding to basis elements in $V_i$, so that the action of
$\gm[k]$ is described by $t\cdot(x_1, \dots, x_r) = (t^{m_1}x_1, \dots,
t^{m_r}x_r)$. Then it is a simple matter to verify that $f$ is given by a
formula of the type
    \[
    f(x_1, \dots, x_r) = \bigl(x_1, x_2 + f_2(x_1), x_3+ f_3(x_1, x_2),
    \dots, x_r + f_r(x_1, \dots, x_{r-1})\bigr)
    \]
and that every polynomial map of this form is an isomorphism.
\end{proof}
\noqed\end{proof}

Now let us show that Theorem~\ref{thm:BBvariant} implies
Proposition~\ref{prop:describesufflimits}. First of all, Theorem
\refall{thm:BBvariant}{2} and a standard argument with the localization
sequence imply that the pullback map $\K(Y_i)
\otimes_\mathbb{Z} \mathrm{R}G = \K(Y_i,G)
\to K(X_i,G)$ is an isomorphism.

Now, let us order the strata $X_1$, \dots,~$X_r$ by decreasing dimension,
and let us set $U_i = X_1 \cup \ldots \cup X_i$. Clearly $X_i$ is closed in
$U_i$. We claim that $X_i$ is K-rigid in $U_i$ for all $i$. In fact, it is
enough to show that the restriction of $\lambda_{-1}(\N_{X_i}X)$ to
$Y_i$ is not a zero-divisor, and this follows from
Lemma~\ref{lem:criterionK-rigid} and Theorem \refall{thm:BBvariant}{3}.

Then by Proposition \refall{prop:K-rigid->}{1} we have an exact sequence
    \[
    0 \longrightarrow \K(X_i,G) \longrightarrow \K(U_i, G) \longrightarrow
    \K(U_{i-1}, G)
    \longrightarrow 0;
    \]
so each $\K(U_i, G)$ is obtained by finitely many successive
extensions of $\R G$-modules of the form $N \otimes_\mathbb{Z}\R G$, and
$X = U_r$. This concludes the proof of
Proposition~\ref{prop:describesufflimits}, and of
Theorem~\ref{thm:->suff-deep}.

\subsection{Comparison with ordinary K-theory for torus actions with
enough limits}

Assume that $T$ is a totally split torus over a perfect field $k$, acting
on a separated scheme $X$ of finite type over $k$. We write $T$ instead of
$G$ for conformity with the standard notation.

The following is a consequence of
Proposition~\ref{prop:describesufflimits}.

\begin{corollary}\label{cor:vanishtor}
If $X$ is smooth and the action has enough limits, we have
    \[
    \tor_p^{\R T}\bigl(\K(X, T), \mathbb{Z}\bigr) = 0
    \text{ for all $p > 0$.}
    \]
\end{corollary}

The interest of this comes from the following result of Merkurjev.

\begin{theorem}[{\cite[Theorem~4.3]{merk}}]\label{thm:merk}
There is a homology spectral sequence
    \[
    E^2_{pq} = \tor_p^{\R T}\bigl(\mathbb{Z}, \K[q](X, T)\bigr)
    \Longrightarrow \K[p+q](X)
    \]
such that the edge homomorphisms
\[
\mathbb{Z}\otimes_{\R T} \K[q](X,T) \longrightarrow \K[q](X)
\]
are induced by the forgetful homomorphism $\K(X,T) \to \K(X)$.

In particular the ring homomorphism $\mathbb{Z}\otimes_{\R T} \K[0](X,T)
\to\K[0](X)$ is an isomorphism.

Furthemore, if $X$ is smooth and projective  we have $E^2_{pq} = 0$ for all
$p>0$, so the homomorphism $\mathbb{Z}\otimes_{\R T} \K(X,T) \to
\K(X)$ is an isomorphism.
\end{theorem}

More generally, Merkurjev produces his spectral sequence for actions of
reductive groups whose fundamental group is torsion-free.

From
Corollary~\ref{cor:vanishtor} we get the following extension of Merkurjev's
degeneracy result.

\begin{theorem}\label{thm:enough->degenerates}
Suppose that  $T$ is a totally split torus over a
perfect field $k$, acting with enough limits on a smooth scheme separated
and of finite type over $k$. Then the forgertful homomorphisms $\K(X, T)
\to \K(X)$ induces an isomorphism
    \[
    \mathbb{Z}\otimes_{\R T} \K(X, T) \larrowsim \K(X).
    \]
\end{theorem}

\section{The K-theory of smooth toric varieties}\label{sec:toric}

Our reference for the theory of toric varieties will be \cite{fultontoric}.

In this section we take $T$ to be a totally split torus over a fixed
field $k$,
    \[
    N = \hom(\gm[k], T) = \widehat T^\vee
    \]
its lattice of one-parameter subgroups, $\Delta$ a fan in
$N \otimes\mathbb{R}$, $X = X(\Delta)$ the associated toric variety. We
will always assume that $X$ is smooth; this is equivalent to saying that
every cone in $\Delta$ is generated by a subset of a basis of $N$.

We will give two different descriptions of the equivariant K-theory ring
of $X$, one as a subring of a product of representation rings, and the
second by generators and relations, analogously to what has been done
for equivariant cohomology in \cite{bdcp}.

\subsection{The equivariant K-theory ring as a subring of a product of
rings of representations}

There is one orbit $O_\sigma$ of $T$ on $X$ for each cone $\sigma\in
\Delta$, containing a canonical rational point $x_\sigma \in O_\sigma(k)$.
The dimension of $O_\sigma$ is the codimension $\codim \sigma \eqdef
\dim T - \dim \sigma$, and the stabilizer of any of its geometric points is
the subtorus $T_\sigma \subseteq T$ whose group of one-parameter subgroups is
precisely the subgroup $\generate \sigma = \sigma + (-\sigma)\subseteq N$;
the dimension of $T_\sigma$ is equal to the dimension of $\sigma$ (see
\cite[3.1]{fultontoric}).  Hence $X_s$ is the disjoint union of the orbits
$O_\sigma = T/T_\sigma$ with $\dim\sigma = s$.

Given a cone $\sigma$ in $N\otimes \mathbb{R}$, we denote by
$\partial\sigma$ the union of all of its faces of codimension~1.

Since
    \[
    \K(O_\sigma, T) = \K(T/T_\sigma, T) = \K(\spec k, T_\sigma) = \K(k)
    \otimes \R T_\sigma
    \]
we have that
    \[
    \prod_{s} \K(X_s, T) = \prod_{\sigma \in \Delta} \K(O_\sigma, T)\
    = \prod_{\sigma \in \Delta} \K(k) \otimes \R T_\sigma.
    \]

\begin{lemma}
Fix an positive integer $s$. Then there is a canonical isomorphism
    \[
    \N_{s, s-1} \simeq \coprod_{\substack{\sigma\in \Delta\\
    \dim \sigma = s}} \coprod_{\tau \in \partial\sigma}
    O_\tau.
    \]

Furthermore, for each pair $\sigma$, $\tau$ such that $\sigma$ has
dimension $s$, $\tau$ has dimension $s-1$, and $\tau$ is a face of
$\sigma$, the composition of the specialization homomorphism
    \begin{align*}
    \sp_{X,s}^{s-1} \colon {}&\K(X_{s-1},T) =
    \prod_{\substack{\tau \in \Delta\\ \dim \tau = {s-1}}}\K(O_\tau,T)\\
    &\qquad\longrightarrow \K(\N_{s, s-1},T) = \prod_{\substack{\sigma\in
    \Delta\\
    \dim \sigma = s}}\prod_{\tau \in \partial\sigma}\K(O_\tau,T),
    \end{align*}
with the projection
    \[
    \pr_{\sigma, \tau}\colon \prod_{\substack{\sigma\in
    \Delta\\ \dim \sigma = s}}
    \prod_{\tau \in \partial\sigma}\K(O_\tau,T)
    \longrightarrow \K(O_\tau,T)
    \]
is the projection
    \[
    \pr_\tau \colon \prod_{\substack{\tau \in \Delta\\
    \dim \tau = {s-1}}}\K(O_\tau,T) \longrightarrow
    \K(O_\tau,T).
    \]
\end{lemma}

\begin{proof}
We us the same notation as in \cite{fultontoric}; in particular, for each
cone $\sigma$ of the fan of $X$, we denote by $U_\sigma$ the corresponding
affine open subscheme of $X$.

First of all, assume that the fan
$\Delta$ consists of all the faces of an
$s$-dimensional cone $\sigma$. Call $B$ a part of a basis of $N$ that
spans $\sigma$: we have an action of $T_\sigma$ on the $k$-vector space
$V_\sigma$ generated by $B$, by letting each 1-parameter subgroup
$\mathbb{G}_\mathrm{m} \to G$ in $B$ act by multiplication on the
corresponding line in
$V_\sigma$, and an equivariant embedding
$T_\sigma \subseteq V_\sigma$. Then $X =U_\sigma$ is $T$-equivariantly
isomorphic to the $T$-equivariant vector bundle
    \[
    T \times^{T_\sigma} V_\sigma = (T \times V_\sigma)/T_\sigma
    \longrightarrow O_\sigma = T/T_\sigma
    \]
in such a way that the zero section corresponds to
$O_\sigma \subseteq U_\sigma$. Since $X_s = O_\sigma$ and $U_\sigma$ is a
vector bundle over $O_\sigma$, we get a canonical isomorphism $U_\sigma
\simeq \N_s$, and from this a canonical isomorphism
    \[
    \N_{s, s-1} \simeq (U_\sigma)_{s-1} = \coprod_{\tau\in
    \partial\sigma}O_\tau.
    \]

It follows that the deformation to the normal bundle of
$O_\sigma$ in $U_\sigma$ is also isomorphic to the product $U_\sigma
\times_k \mathbb{P}^1$, and from this we get the second part of the
statement.

In the general case we have $X_s = \coprod_{\dim\sigma = s}O_\sigma$, and
if $\sigma$ is a cone of dimension $s$ in $\Delta$, the intersection of
$X_s$ with $U_\sigma$ is precisely $O_\sigma$. From this we get the first
part of the statement in general.

The second part follows by applying the compatibility of specializations
to the morphism of deformation to the normal bundle induced by the
equivariant morphism $\coprod_{\dim\sigma = s}U_\sigma \to X_s$.
\end{proof}

Using this lemma together with Theorem~\ref{thm:maintheorem} we get that
$\K(X, T)$ is the subring of
    \[
    \prod_{\sigma\in \Delta} \K(O_\sigma,T) =
    \prod_{\sigma\in \Delta} \K(k) \otimes \R T_\sigma
    \]
consisting of elements $(a_\sigma)$ with the property that the restriction
of $a_\sigma \in \K(k)\otimes \R T_\sigma$ to $ \K(k)\otimes
\R T_\tau$ coincides with $a_\tau \in \K(k)\otimes \R T_\tau$ every time
$\tau$ is a face of codimension~1 in $\sigma$. Since every face of a cone
is contained in a face of codimension~1, this can also be described as the
subring of $\prod_{\sigma} \K(k) \otimes \R T_\sigma$ consisting of
elements $(a_\sigma)$ with the property that the restriction of $a_\sigma
\in \K(k)\otimes \R T_\sigma$ to $ \K(k)\otimes
\R T_\tau$ coincides with $a_\tau \in \K(k)\otimes \R T_\tau$ every time
$\tau$ is a face of $\sigma$. But every cone in $\Delta$ is contained in a
maximal cone in $\Delta$, so we get the following description of the
equivariant K-theory of a smooth toric variety.

\begin{theorem}\label{thm:describetoric}
If $X(\Delta)$ is a smooth toric variety associated with a fan $\Delta$
in $N \otimes\mathbb{R}$, there is an injective homomorphism of
$\R T$-algebras
    \[
    \K\bigl(X(\Delta),T\bigr) \into \prod_{\sigma \in \Delta_{\max}}\K(k)
    \otimes \R T_\sigma,
    \]
where $\Delta_{\max}$ is the set of maximal cones in $\Delta$.

An element $(a_\sigma) \in \prod_{\sigma}\K(k) \otimes \R T_\sigma$ is in
the image of this homomorphism if and only if for any two maximal cones
$\sigma_1$ and
$\sigma_2$, the restrictions of $a_{\sigma_1}$ and $a_{\sigma_2}$ to $\K(k)
\otimes \R T_{\sigma_1\cap \sigma_2}$ coincide.
\end{theorem}

This description of the ring $\K(X, T)$ is analogous to the description of
its equivariant cohomology in \cite{bdcp}, and of it equivariant Chow ring in
\cite[Theorem~5.4]{brion1}.

\subsection{The multiplicative Stanley--Reisner presentation}
\label{subsec:multSR}

 From the description above it is easy to get a presentation of the
equivariant K-theory ring of the smooth toric variety $X(\Delta)$,
analogous to the Stanley--Reisner presentation of its equivariant
cohomology ring obtained in \cite{bdcp}. Denote by $\Delta_1$ the
subset of $\Delta$ consisting of 1-dimensional cones. We will use the
following notation: if
$\sigma \in \Delta_{\max}$, call
$N_\sigma \subseteq N$ the group of 1-parameter subgroups of $T_\sigma$,
so that $\widehat T_\sigma = (N_\sigma)^\vee$. We will use
multiplicative notation for $\widehat T_\sigma$. Furthermore, for any
$\rho \in \Delta_1$ we call $v_\rho \in N$ the generator for the
monoid  $\rho \cap N$.

For each $\rho \in \Delta_1$ we define an element $u_\rho$ of the product
$\prod_{\sigma\in \Delta_{\max}} \widehat T_{\sigma}$, as follows. If
$\rho$ is not a face of $\sigma$, we set $(u_\rho)_\sigma = 1$. If $\rho \in
\Delta_{1}$ is a face of $\sigma \in \Delta_{\max}$ and $\{\rho =
\rho_1, \rho_2\dots, \rho_t\}$ is the set of 1-dimensional faces of
$\sigma$, then, since the variety $X(\Delta)$ is smooth, we have that
$v_{\rho_1}$, \dots,~$v_{\rho_t}$ form a basis for $N_\sigma$. If
$v^\vee_{\rho_1}$, \dots,~$v^\vee_{\rho_t}$ is the dual basis in $\widehat
T_\sigma = N_\sigma^\vee$, we set $u_\rho = v^\vee_{\rho_1}$.

Denote by $V_\Delta \subseteq \prod_{\sigma\in \Delta_{\max}} \widehat
T_\sigma$ the subgroup consisting of the elements $(x_\sigma)$ with the
property that for all
$\sigma_1$, $\sigma_2$ in $\Delta_{\max}$ the restrictions of
$x_{\sigma_1} \in \widehat T_{\sigma_1}$ and $x_{\sigma_2} \in \widehat
T_{\sigma_2}$ to $\widehat T_{\sigma_1 \cap \sigma_2}$ coincide. Then we
have the following fact.

\begin{proposition}\label{prop:basis}
The elements $u_\rho$ form a basis of $V_\Delta$.
\end{proposition}

The proof is straightforward.

We have the inclusions
    \[
    \prod_{\sigma \in \Delta_{\max}} \widehat T_\sigma
    \subseteq \prod_{\sigma \in \Delta_{\max}} \R T_\sigma
    \subseteq \prod_{\sigma \in \Delta_{\max}}
    \K(k) \otimes \R T_\sigma;
    \]
because of the description of the ring $\K\bigl(X(\Delta), T\bigr)$ as a
subring of $\prod_{\sigma \in \Delta_{\max}} \K(k) \otimes \R
T_\sigma$ given in Theorem~\ref{thm:describetoric}, we see that we can
consider the $u_\rho$ as elements of $\K\bigl(X(\Delta), T\bigr)$.

There are some obvious relations that the $u_\rho$ satisfy in
    \[
    \K\bigl(X(\Delta), T\bigr) \subseteq \prod_{\sigma \in \Delta_{\max}}
    \K(k) \otimes \R T_\sigma.
    \]
Suppose that $S$ is a subset of $\Delta_1$ not
contained in any maximal cone of $\Delta$. Then for all $\sigma$ in
$\Delta_{\max}$ there will be some $\rho \in S$ such that $(u_\rho)_\sigma
= 1$ in $\widehat T_\sigma$; hence we have the relation
    \[
    \prod_{\rho \in S} (u_\rho - 1) = 0 \quad
    \text{in} \quad \K\bigl(X(\Delta), T\bigr) \subseteq \prod_{\sigma \in
    \Delta_{\max}} \K(k) \otimes \R T_\sigma.
    \]
 From this we get a homomorphism of $\K(k)$-algebras
    \begin{equation} \label{eq:hom}
    \frac{\K(k)\bigl[x_\rho^{\pm 1}\bigr]} {\bigl(\prod_{\rho \in S}
    (x_\rho - 1)\bigr)} \longrightarrow \K\bigl(X(\Delta), T\bigr)
    \end{equation}
where the $x_\rho$ are indeterminates, where $\rho$ varies over
$\Delta_1$, and $S$ over the subsets of $\Delta_1$ whose elements
are not all contained in a maximal cone in $\Delta$, by sending each
$x_\rho$ to $u_\rho$.

\begin{theorem}\label{thm:SR}
Suppose that $X(\Delta)$ is a smooth toric variety associated with a fan
$\Delta$ in $N \otimes\mathbb{R}$. Then the homomorphism~(\ref{eq:hom})
above is an isomorphism.
\end{theorem}

\begin{proof}
First of all, let us show that the $u_\rho$ and their inverses generate
    \[
    \K\bigl(X(\Delta), T\bigr) \subseteq \prod_{\sigma \in \Delta_{\max}}
    \K(k) \otimes\R T_\sigma.
    \]
Set $\Delta_{\max} = \{\sigma_1, \dots,
\sigma_r\}$, and let $\alpha$ be an element of
$\K\bigl(X(\Delta), T\bigr)$; we want to show that $\alpha$ can be
expressed as a Laurent polynomial in the $x_\rho$ evaluated in the
$u_\rho$. The ring $\K(k)
\otimes \R T_{\sigma_1}$ is a ring of Laurent polynomials in the images of
the $u_\rho$ with $\rho \subseteq \sigma_1$, so we can find a Laurent
polynomial $p_1(x_\rho)$, in which only the $x_\rho$ with $\rho \subseteq
\sigma_1$ appear, such that the image of $p_1(u_\rho)$ in $\K(k)
\otimes \R T_{\sigma_1}$ equals the image of $\alpha$ in the same ring. By
subtracting $p_1(u_\rho)$, we see that we may assume that the projection of
$\alpha$ into $\K(k) \otimes \R T_{\sigma_1}$ is zero.

Now, let us repeat the procedure for the maximal cone $\sigma_2$: find a
polynomial $p_1(x_\rho)$, in which only the $x_\rho$ with $\rho \subseteq
\sigma_2$ appear, such that the image of $p_2(u_\rho)$ in $\K(k)
\otimes \R T_{\sigma_2}$ equals the image of $\alpha$ in the same ring.
The key point is that the restriction of $\alpha$ to $\K(k)
\otimes \R T_{\sigma_1 \cap \sigma_2}$ is zero, so in fact $p_2(x_\rho)$
can only contain the variables $x_\rho$ with $\rho$ not in $\sigma_1$.
Hence the restriction of $p_2(x_\rho)$ to $\K(k) \otimes \R T_{\sigma_1}$
is also zero, and after having subtracted $p_2(x_\rho)$ from $\alpha$ we
may assume that the restriction of $\alpha$ to both $\K(k) \otimes \R
T_{\sigma_1}$ and $\K(k) \otimes \R T_{\sigma_2}$ is zero. We can continue
this process for the remaining cones $\sigma_3$, \dots,~$\sigma_r$; at the
end all the projections will be zero, and therefore $\alpha$ will be zero
too.

Now we have to show that the kernel of the homomorphism
    \begin{equation}\label{eq:homo2}
    \K(k)\bigl[x_\rho^{\pm 1}\bigr] \longrightarrow \prod_{\sigma
    \in \Delta_{\max}} \K(k)\otimes \R T_\sigma
    \end{equation}
sending each $x_\rho$ to $u_\rho$ equals the ideal $\bigl(\prod_{\rho \in
S} (x_\rho - 1)\bigr)$, where $S$ varies over all subsets of $\Delta_1$
not contained in any maximal cone. The kernel of the
projection
    \[
    \K(k)\bigl[x_\rho^{\pm 1}\bigr] \longrightarrow
    \K(k)\otimes \R T_\sigma
    \]
equals the ideal $I_\sigma$ generated by the $x_\rho - 1$, where $\rho$
varies over the set of 1-dimensional cones in $\Delta_1$ not contained in
$\sigma$, hence the kernel of the homomorphism~\ref{eq:homo2} is the
intersection of the $I_\sigma$. The result is then a consequence of the
following lemma.

\begin{lemma}
Let $R$ be a (not necessarily commutative) ring, $\{x_\rho\}_{\rho \in E}$
a finite set of central indeterminates; consider the ring of Laurent
polynomials $R\bigl[x_\rho^{\pm 1} \bigr]$. Let $A_1$,
\dots,~$A_r$ be subsets of $E$; for each $j = 1$, \dots,~$r$ call $I_j$
the ideal of $R\bigl[x_\rho^{\pm 1} \bigr]$ generated by the elements
$x_\rho - 1$ with $\rho \in A_j$.

Then the intersection $I_1 \cap \dots \cap I_r$ is the ideal of
$R\bigl[x_\rho^{\pm 1} \bigr]$ generated by the elements $\prod_{\rho \in
S}(x_\rho - 1)$, where $S$ varies over all
subsets of $E$ that meet each $A_j$.
\end{lemma}

\begin{remark}
When each $A_j$ contains a single element this is a particular case of
Lemma~\ref{lem:int=prod}. The obvious common generalization should also hold.
\end{remark}

\begin{proof}

We proceed by induction on $r$; the case $r=1$ is clear. In general, take
$p \in I_1 \cap \dots \cap I_r \subseteq I_2 \cap \dots \cap I_r$; by
induction hypothesis, we can write
    \[
    p = \sum_S \biggl(\prod_{\rho \in S}\
    (x_\rho - 1) \biggr) q_S,
    \]
where $S$ varies over all subsets of $E$ whose
intersection with each of the $A_2$, \dots,~$A_r$ is not empty. We can
split the sum as
    \[
    p = \sum_{S\cap A_1 \neq \emptyset} \biggl(\prod_{\rho \in S}
    (x_\rho - 1) \biggr) q_S + \sum_{S\cap A_1 =
    \emptyset}\biggl(\prod_{\rho \in S}(x_\rho - 1) \biggr) q_S;
    \]
the first summand is in $I_1 \cap \dots \cap I_r$ and is of the desired
form, so we may subtract it from $p$ and suppose that $p$ is of the type
    \[
    p = \sum_{S\cap A_1 =
    \emptyset}\biggl(\prod_{\rho \in S}(x_\rho - 1) \biggr) q_S.
    \]
Now, consider the ring $R\bigl[x_\rho^{\pm 1} \bigr]_{\rho \notin A_1}$ of
Laurent polynomials not involving the variables in $A_1$; it is a subring
of $R\bigl[x_\rho^{\pm 1} \bigr]$, and there is also a retraction
    \[
    \pi \colon R\bigl[x_\rho^{\pm 1} \bigr] \longrightarrow
    R\bigl[x_\rho^{\pm 1}\bigr]_{\rho \notin A_1},
    \]
sending each $x_\rho$ with $\rho \in A_1$ to $1$, whose kernel is
precisely $I_1$. The elements $\prod_{\rho \in S}(x_\rho - 1)$ are in
$R\bigl[x_\rho^{\pm 1}\bigr]_{\rho \notin A_1}$, and
    \[
    \pi p = \sum_{S\cap A_1 =
    \emptyset}\biggl(\prod_{\rho \in S}(x_\rho - 1) \biggr) \pi q_S = 0
    \in R\bigl[x_\rho^{\pm 1}\bigr]_{\rho \notin A_1}
    \]
so we can write
    \[
    p = \sum_{S\cap A_1 =
    \emptyset}\biggl(\prod_{\rho \in S}(x_\rho - 1) \biggr) (q_S - \pi q_S).
    \]
Then we write each $q_S - \pi q_S \in I_1$ as a linear combination of the
polynomials $x_\rho - 1$ with $\rho \in A_1$; this concludes the proofs
of the lemma and of the theorem.

\end{proof}

\noqed\end{proof}

\subsection{Ordinary K-theory of smooth toric varieties}

 From Merkurjev's theorem (\ref{thm:merk}) we get that $\K[0](X) =
\mathbb{Z}\otimes_{\R T} \K[0](X, T)$. If the Merkurjev spectral
sequence degenerates then we also have  $\K(X) = \mathbb{Z}\otimes_{\R T}
\K(X, T)$; this gives a way to compute the whole K-theory ring of $X$.

In general, the spectral sequence will not degenerate, and the ring
$\K(X)$ tends to be rather complicated (for example, when $X = T$). When the
toric variety has enough limits we can apply
Theorem~\ref{thm:enough->degenerates}. Using the description of closures of
orbits given in
\cite[3.1]{fultontoric}, one shows that given a point $x\in
X\bigl(\,\overline k\,\bigr)$ lying in an orbit $O_\tau$, and a
one-parameter subgroup corresponding to an element $v \in N$, then the
point has a limit under the one parameter subgroup if and only if $v$ lies
in the subset
    \[
    \bigcup_{\sigma \in \star \tau} (\sigma + \generate{\tau})
    \subseteq N \otimes \mathbb{R},
    \]
where $\generate{\tau}$ denotes the subvector space
$\tau + (-\tau) \subseteq N \otimes \mathbb{R}$, and $\star \tau$ is the
set of cones in $\Delta$ containing $\tau$ as a face. From this we obtain
the followign.

\begin{proposition}\label{charact}
$X$ has enough limits if and only if the subset
    \[
    \bigcap_{\tau\in \Delta}\,
    \bigcup_{\sigma \in \star \tau} (\sigma + \generate{\tau})
    \subseteq N \otimes \mathbb{R}
    \]
has nonempty interior.
\end{proposition}

\begin{remark}\label{rmk:describecombcomplete} If $X = X(\Delta)$ is a
smooth toric variety with enough limits, the K-theory ring of $X$
can be describe in a slightly more efficient fashion: there is an
injective homomorphism of $\R T$-algebras
    \[
    \K(X,T) \into \prod_{\substack{\sigma \in \Delta\\
    \dim \sigma = \dim T}}\K(k) \otimes \R T,
    \]
and an element $(a_\sigma) \in \prod_{\sigma}\K(k) \otimes \R T$ is
in the image of this homomorphism if and only if for any two
\emph{adjacent} maximal cones
$\sigma_1$ and $\sigma_2$, the restrictions of $a_{\sigma_1}$ and
$a_{\sigma_2}$ to $\K(k)\otimes \R T_{\sigma_1\cap \sigma_2}$ coincide.
\end{remark}

\begin{theorem}\label{thm:combcomplete->}
If $X$ is a smooth smooth toric variety with enough limits, then
$\K[0](X,T)$ is a projective module over $\R T$ of rank equal to the
number of maximal cones in its fan; furthermore the natural ring
homomorphism $\K(k) \otimes \K[0](X,T) \to \K(X,T)$ is an isomorphism.
\end{theorem}

In particular we have
    \[
    \tor_p^{\R T}\bigl(\K(X, T), \mathbb{Z}\bigr) = 0
    \text{ for all $p > 0$:}
    \]
so from Merkurjev's theorem (\ref{thm:merk}) we get the following.

\begin{corollary}
Let $X$ be a smooth toric variety with enough limits.

\begin{enumeratei}

\item The natural
homomorphism of rings
    \[
    \mathbb{Z}\otimes_{\R T} \K(X, T) \longrightarrow \K(X)
    \]
is an isomorphism of $\R T$-algebras.

\item $\K[0](X)$ is a free
abelian group of rank equal to the number of maximal cones in $\Delta$.

\item The natural homomorphism $\K(k) \otimes \K[0](X) \to \K(X)$ is an
isomorphism.

\end{enumeratei}
\end{corollary}

\begin{remark}
Suppose that the base field is the field $\mathbb{C}$ of complex numbers. The
Merkurjev spectral sequence is an analogue of the Eilenberg--Moore spectral
sequence (\cite{em})
    \[
    E_2^{p,q} = \tor_{p,q}^{H^*(BT, \mathbb{Z})}\bigl(H^*_T(X,\mathbb{Z}),
    \mathbb{Z}\bigr) \Longrightarrow \H^{p+q}(X, \mathbb{Z}).
    \]
Then \cite{bbfk} contains a description of the fans of the simplicial toric
varieties for which this spectral sequence degenerates after tensoring with
$\mathbb{Q}$. Presumably there should be a similar description for the case
considered here.
\end{remark}

\section{The refined decomposition theorem}\label{sec:decomposition}

The main result of \cite{vevi} shows that if $G$ is an algebraic group
acting with finite stabilizers on a noetherian regular algebraic space
$X$ over a field, the equivariant K-theory ring of $X$, after inverting
certain primes, splits as a direct product of rings related with the
K-theory of certain fixed points subsets. For actions of diagonalizable
groups it is not hard to extend this decomposition to the case when the
stabilizers have constant dimension.

So, in the general case when we do not assume anything about the dimension
of the stabilizers, this theorem gives a description of the K-theory of
each stratum $X_s$; it should clearly be possible to mix this with
Theorem~\ref{thm:maintheorem} to give a result that expresses $\K(X,G)$,
after inverting certain primes, as a fibered product. This is carried out
in this section.

The material in this section is organized as follows. We first extend the
main result of \cite{vevi}  over an arbitrary noetherian separated base
scheme $S$ (Theorem~\ref{thm:refined0}) by giving a decomposition theorem in
the case of an action with finite stabilizers of a diagonalizable
group scheme $G$ of finite type over $S$ on a noetherian regular separated
algebraic space $X$ over $S$. Next, we deduce from this a decomposition
theorem in the case where the action of $G$ on $X$ has stabilizers
of fixed constant dimension (Theorem~\ref{thm:refinedconstant}). Finally,
we combine the analysis carried out in Section 4, and culminating with
Theorem~\ref{thm:maintheorem}, together with
Theorem~\ref{thm:refinedconstant} to prove a general decomposition theorem
(Theorem~\ref{thm:refineddecomposition}) where no restriction is
imposed on the stabilizers.

\subsection{Actions with finite stabilizers}

Here we recall the main result of \cite{vevi} for actions of
diagonalizable groups, extending it over any noetherian separated 
base scheme $S$.

Suppose that $G$ is a diagonalizable group scheme of finite type acting
with finite stabilizers on a noetherian regular separated algebraic space
over a noetherian separated scheme $S$. A diagonalizable group scheme of
finite type
$\sigma$ over $S$ is called
\emph{dual cyclic} if its Cartier dual is finite cyclic, that is, if
$\sigma$ is isomorphic to a group scheme of the form
$\boldsymbol{\mu}_{n,S}$ for some positive integer $n$.

A subgroup scheme $\sigma\subseteq G$ is called \emph{essential} if it is
dual cyclic, and $X^\sigma \neq \emptyset$.

There are only finitely many essential subgroups of $G$; we will fix a
positive integer $N$ which is divisible by the least common multiple of
their orders.

Suppose that $\sigma$ is a dual cyclic group of order $n$. The ring of
representations $\R \sigma$ is of the form $\mathbb{Z}[t]/(t^n - 1)$,
where $t$ corresponds to a generator of the group of characters
$\widehat\sigma$. Denote by $\Rtilde \sigma$ the quotient of $\R \sigma$
corresponding to the quotient
    \[
    \R \sigma = \mathbb{Z}[t]/(t^n - 1) \twoheadrightarrow
    \mathbb{Z}[t]/\bigl(\Phi_n(t)\bigr),
    \]
where $\Phi_n$ is the $n\th$ cyclotomic polynomial. This quotient $\Rtilde
\sigma$ is independent of the choice of a generator for $\widehat\sigma$.
We have a canonical homomorphism $\R G \to \R \sigma \twoheadrightarrow
\Rtilde \sigma$ induced by the embedding $\sigma \subseteq G$.

We also define a multiplicative system
    \[
    \mathrm{S}_\sigma \subseteq \R G
    \]
as follows: an element of $\R G$ is in $\mathrm{S}_\sigma$ if its
image in $\Rtilde \sigma$ is a power of $N$.

For any $\R G$-module $M$, we define the \emph{$\sigma$-localization}
$M_\sigma$ of $M$ to be $\mathrm{S}_\sigma^{-1} M$.  Consider the
$\sigma$-localization
    \[
    \K(X^\sigma,G)_\sigma = \mathrm{S}_\sigma^{-1} \K(X^\sigma,G)
    \]
of the  $\R G$-algebra
$\K(X^\sigma,G)$. The tensor product
$\K(X^\sigma,G)_\sigma \otimes \mathbb{Q}$ is the localization
    \[
    \bigl(\K(X^\sigma,G)\otimes \mathbb{Q}\bigr)_{\mathfrak{m}_\sigma}
    \]
of the $\R \sigma$-algebra $\K(X^\sigma, G)\otimes\mathbb{Q}$ at the
maximal ideal
    \[
    \mathfrak{m}_\sigma  = \ker(\R \sigma\otimes\mathbb{Q}
    \twoheadrightarrow \Rtilde \sigma \otimes \mathbb{Q}).
    \]
We are particularly interested in the $\sigma$-localization when $\sigma$
is the trivial subgroup of $G$; in this case we denote it by
$\K(X,G)\geom$, and call it \emph{the geometric equivariant K-theory of
$X$}. The localization homomorphism
    \[
    \K(X,G)\otimes \mathbb{Z}[1/N] \longrightarrow \K(X,G)\geom
    \]
is surjective, and its kernel can be described as follows. Consider the
kernel $\mathfrak{p} = \ker \rk$ of the localized rank homomorphism
    \[
    \rk \colon \K(X,G)\otimes\mathbb{Z}[1/N] \longrightarrow
    \mathbb{Z}[1/N];
    \]
then the power $\mathfrak{p}^k$ is independent of $k$ if $k$ is large, and
this power coincides with the kernel of the localization homomorphism.

For each essential subgroup $\sigma\subseteq G$, consider the
compositions
    \[
    \loc_\sigma \colon \K(X,G)\otimes \mathbb{Z}[1/N] \longrightarrow
    \K(X^\sigma,G)\otimes\mathbb{Z}[1/N] \longrightarrow
    \K(X^\sigma,G)_\sigma,
    \]
where the first arrow is a restriction homomorphism, and the second one is
the localization.

There is also a homomorphism of $\R G$-algebras $\K(X^\sigma, G) \to
\K(X^\sigma, G)\geom \otimes \Rtilde \sigma$, defined as the composition
    \[
    \K(X^\sigma, G) \longrightarrow \K(X^\sigma, G\times\sigma)
    \simeq \K(X^\sigma, G)\otimes \R \sigma \longrightarrow
    \K(X^\sigma, G)\geom \otimes\Rtilde \sigma,
    \]
where the first morphism is induced by the multiplication $G \times \sigma
\to G$, the second one is a natural isomorphism coming from the fact that
$\sigma$ acts trivially on $X^\sigma$ (\cite[Lemma~2.7]{vevi}), and the
third one is obtained from the localization homomorphism $\K(X^\sigma, G)
\to \K(X^\sigma, G)\geom$ and the projection $\R G \to \Rtilde \sigma$.
Then the homomorphism $\K(X^\sigma, G) \to
\K(X^\sigma, G)\geom \otimes \Rtilde \sigma$ factors through 
$\K(X^\sigma, G)_\sigma$
(\cite[Lemma~2.8]{vevi}), inducing a homomorphism
    \[
    \theta_\sigma\colon \K(X^\sigma, G)_\sigma \longrightarrow\K(X^\sigma,
    G)\geom
    \otimes\Rtilde \sigma.
    \]

\begin{theorem}\call{thm:refined0}\hfil

\begin{enumeratei}

\itemref{1} There are finitely many essential subgroup schemes in $G$, and
the homomorphism
    \[
    \prod_\sigma \loc_\sigma\colon \K(X,G)\otimes \mathbb{Z}[1/N]
    \longrightarrow \prod_\sigma \K(X^\sigma,G)_\sigma,
    \]
where the product runs over all the essential subgroup schemes of
$G$, is an isomorphism.

\itemref{2}  The homomorphism
    \[
    \theta_\sigma\colon \K(X^\sigma, G)_\sigma \longrightarrow
    \K(X^\sigma, G)\geom
    \otimes\Rtilde \sigma
    \]
is an isomorphism of $\R G$-algebras.

\end{enumeratei}
\end{theorem}

\begin{proof}
If the base scheme $S$ is the spectrum of a field, this is a particular
case of the main theorem of \cite{vevi}.

If $G$ is a torus, the proof of this statement given in \cite{vevi} goes
through without changes, because it only relies on Thomason's generic slice
theorem for torus actions (\cite[Proposition~4.10]{th2}).

In the general case, choose an
embedding $G\into T$ into some totally split torus $T$ over $S$, and
consider the quotient space
    \[
    Y \eqdef X\times ^GT = (X \times T)/G
    \]
by the customary diagonal action of $G$; this exists as an algebraic
space thanks to a result of Artin (\cite[Corollaire~10.4]{lmb}). The same
argument as in the beginning of Section~5.1 of \cite{vevi} shows that
$Y$ is separated.

Now observe that if $\sigma \subseteq T$ is an
essential subgroup relative to the action of $T$ on $Y$, we
have
$Y^\sigma =\emptyset $ unless $\sigma \subseteq G$ is an essential
subgroup; hence the least common multiples of the orders of all essential
subgroups are the same for the action of $G$ on $X$ and the action of $T$ on
$Y$.

Also, if $\sigma \subseteq G$ is an essential
subgroup we have $Y^\sigma =X^\sigma
\times ^GT$,
and therefore, by Morita equivalence (\cite[Proposition~6.2]{th5}), we get
an isomorphism
    \[
    \K(Y^\sigma ,T)\simeq \K(X^\sigma ,G)
    \]
which is an isomorphism of $\R T$-algebras, if we view $\K(X^\sigma ,G)$
as an $\R T$-algebra via the restriction homomorphism $\R T \to \R G$.

Moreover, $\mathrm{S}_\sigma
^T\subseteq \mathrm{R}(T)$ is exactly the preimage of $\mathrm{S}_\sigma
^G\subseteq \mathrm{R}(G)$ under the natural surjection
$\R T\to \R G$; therefore we have compatible $\sigma $-localized Morita
isomorphisms
    \[
    \left(\mathrm{S}_\sigma ^T\right)^{-1}\K(Y^\sigma,T)\simeq
    \left(\mathrm{S}_\sigma ^G\right)^{-1}\K(X^\sigma,G)
    \]
and (for $\sigma$ equal to the trivial subgroup)
    \[
    \K(Y^\sigma
    ,T)_{\mathrm{geom}}\simeq \K(X^\sigma ,G)_{\mathrm{geom}};
    \]
hence the theorem for the action of $G$ on $X$ follows from the theorem
for the action of $T$ on $Y$.
\end{proof}

\subsection{Actions with stabilizers of constant dimension}

 From the theorem on actions with finite stabilizers we can easily get a
decomposition result when we assume that the stabilizers have constant
dimension. Assume that $G$ is a diagonalizable group scheme of finite type
over $S$, acting on a noetherian regular separated algebraic space $X$ over
$S$ with stabilizers of constant dimension equal to $s$.

\begin{definition}
A diagonalizable subgroup scheme
$\sigma\subseteq G$ is \emph{dual semicyclic} if $\sigma/\sigma_0$ is dual
cyclic, where $\sigma_0$ is the toral component of $\sigma$.

The \emph{order} of a dual cyclic group $\sigma$ is by definition equal the
order of $\sigma/ \sigma_0$.
\end{definition}

Equivalently,
$\sigma\subseteq G$ is dual semicyclic if it is isomorphic to $\gm^r \times
\boldsymbol{\mu}_{n,S}$ for some $r \ge 0$ and $n > 0$.

\begin{definition}
A subgroup scheme $\sigma \subseteq G$
is called \emph{essential} if it is dual semicyclic and $s$-dimensional,
and $X^\sigma \neq \emptyset$.
\end{definition}

There are finitely many subtori $T_j \subseteq G$ of dimension $s$ in
$G$ with $X^{T_j} \neq \emptyset$, and $X$ is the disjoint union of the
$X^{T_j}$. The toral part of an essential subgroup of $G$
coincides with one of the $T_j$; hence there are only finitely many
essential subgroups of $G$.

We fix a positive integer $N$ which is divisible by the least common
multiple of the orders of the essential subgroups of $G$.

For each dual semicyclic subgroup $\sigma \subseteq G$, we define a
multiplicative system
    \[
    \mathrm{S}_\sigma \eqdef \mathrm{S}_{\sigma/
    \sigma_0} \subseteq \R(G/\sigma_0) \subseteq \R G.
    \]
as the set of those elements of $\R(G/\sigma_0)$ whose image in $\Rtilde
(\sigma/\sigma_0)$ is a power of $N$. If $M$ is a module over $\R G$, we
define, as before, the
$\sigma$-localization of $M$ to be
$M_\sigma =  \mathrm{S}_\sigma^{-1}M$.

If $\sigma\subseteq G$ is an essential subgroup, we can
choose a splitting $G \simeq (G/\sigma_0) \times \sigma_0$; according to
\cite[Lemme~5.6]{th4} this splitting induces an isomorphism
    \[
    \K(X^\sigma, G) \simeq \K(X^\sigma, G/\sigma_0) \otimes \R \sigma_0
    \]
and also an isomorphism of $\sigma$-localizations
    \[
    \K(X^\sigma, G)_\sigma
    \simeq \K(X^\sigma, G/\sigma_0)_{\sigma/\sigma_0} \otimes \R\sigma_0.
    \]
Fix one of the $T_j$, and choose a splitting $G \simeq G/T_j \times T_j$.
We have a commutative diagram
    \[
    \xymatrix{
    {}\K(X^{T_j}, G)\otimes \mathbb{Z}[1/N]\ar[r] \ar[d]^{\sim} &
    {}\prod\limits_{\substack{\sigma\text{ essential}\\ \sigma_0 = T_j}}
    \ar[d]^{\sim}
    {}\K(X^\sigma, G)_\sigma\\
    {}\K(X^{T_j}, G/T_j) \otimes \R T_j \otimes \mathbb{Z}[1/N]
    \ar[r]^-{\sim}
    & {}\prod\limits_{\substack{\sigma\text{ essential}\\
    \sigma_0 = T_j}}
    {}\K(X^\sigma, G/\sigma_0)_{\sigma/\sigma_0 } \otimes \R \sigma_0
    }
    \]
where the two columns are isomorphisms induced by the choice of a splitting
$G \simeq G/T_j \times T_j$, and the rows are induced by composing the
restriction homomorphism from $X^{T_j}$ to $X^\sigma$ with the localization
homomorphism. The bottom row is in an isomorphism because of
Theorem \refall{thm:refined0}{1}.

Since the product of the restriction homomorphisms
    \[
    \K(X,G) \longrightarrow \prod_{j} \K(X^{T_j}, G)
    \]
is an isomorphism, we obtain the following generalization of
Theorem~\ref{thm:refined0}.

\begin{theorem}\call{thm:refinedconstant} Suppose that a diagonalizable
group scheme $G$ of finite type over $S$ acts with stabilizers of 
constant dimension
on a noetherian regular separated algebraic space $X$ over $S$.

\begin{enumeratei}

\itemref{1} There are finitely many essential subgroup
schemes in $G$, and the homomorphism
    \[
    \prod_\sigma \loc_\sigma\colon \K(X,G)\otimes \mathbb{Z}[1/N]
    \longrightarrow \prod_\sigma \K(X^\sigma,G)_\sigma,
    \]
where the product runs over all the essential subgroup schemes of
$G$, is an isomorphism.

\itemref{2}  For any essential subgroup
scheme $\sigma\subseteq G$, a choice of a splitting $G \simeq
(G/\sigma_0)
\times \sigma_0$ gives an isomorphism
    \[
    \K(X^\sigma, G)_\sigma \longrightarrow
    \K(X^\sigma, G/\sigma_0)\geom\otimes\Rtilde \sigma \otimes \R \sigma_0.
    \]

\end{enumeratei}
\end{theorem}

\begin{remark} If $s = 0$, then $\sigma_0 = 1$ for each essential subgroup
$\sigma\subseteq G$, so there is a unique splitting $G \simeq (G/\sigma_0)
\times \sigma_0$, and the isomorphism in \refpart{thm:refinedconstant}{2}
is canonical.
\end{remark}

\subsection{More specializations}

For the refined decomposition theorem we need more specialization
homomorphisms. Let a diagonalizable group scheme $G$ of finite type 
over $S$ act
on a noetherian regular separated algebrac space $X$ over $S$, with 
no restriction
on the dimensions of the stabilizers.

\begin{notation} Given a diagonalizable subgroup scheme $\sigma\subseteq
G$, we set
    \[
    X^{(\sigma)} = X^{\sigma} \cap X_{\le\dim \sigma}.
    \]
\end{notation}

Equivalently, $X^{(\sigma)} = (X_{\dim \sigma})^{\sigma}$.

Obviously $X^{(\sigma)}$ is  a locally closed regular subspace of $X$.

\begin{definition}
Let $\sigma$ and $\tau$ be two diagonalizable subgroup schemes of $G$. We
say that $\tau$ is \emph{subordinate to $\sigma$}, and we write
$\tau\prec \sigma$, if $\tau$ is contained in $\sigma$, and the induced
morphism $\tau \to \sigma/\sigma_0$ is surjective.
\end{definition}

Suppose that $\sigma$ and $\tau$ are diagonalizable subgroup schemes of
$G$ of dimension $s$ and $t$ respectively, and that $\tau$ is subordinate
to $\sigma$. Consider the deformation to the normal cone $\mathrm{M}_s
\to \mathbb{P}^1_S$ of $X_s$ in $X_{\le s}$, considered in
Subsection~\ref{subsec:specializationsnormal}. By
Proposition~\ref{prop:restrictregularatinfinity}, the restriction
$\mathrm{M}_s^{(\tau)} \to \mathbb{P}^1_S$ is regular at infinity, so we
can define a specialization homomorphism
    \[
    \K(X^{(\tau)}, G) \longrightarrow \K(\N_s^{(\tau)}, G).
    \]
Denote by $\N_{\sigma}$ the restriction of $\N_s$ to $X^{(\sigma)}$. We
define the specialization homomorphism
    \[
    \sp_{X, \sigma}^\tau \colon \K(X^{(\tau)}, G) \longrightarrow
    \K(\N_\sigma^{(\tau)}, G)
    \]
as the composition of the homomorphism $\K(X^{(\tau)},G) \to
\K(\N_s^{(\tau)}, G)$ above with the restriction homomorphism
$\K(\N_s^{(\tau)}, G) \to \K(\N_\sigma^{(\tau)}, G)$.

We also denote by
     \[
    \sp_{X, \sigma}^\tau \colon \K(X^{(\tau)}, G)_\tau \longrightarrow
    \K(\N_\sigma^{(\tau)}, G)_\tau
    \]
the $\tau$-localization of this specialization homomorphism.

\begin{remark}
Since $\tau$ is subordinate to $\sigma$, it is easy to see that
$\N_\sigma^{(\tau)}$ is a union of connected components of $\N_s^{(\tau)}$.
\end{remark}

\subsection{The general case}

The hypotheses are the same as in the previous subsection: $G$ is a
diagonalizable group scheme of finite type over $S$, acting on a
noetherian regular separated algebraic space $X$ over $S$.

\begin{definition}
An \emph{essential subgroup of $G$} is a dual semicyclic subgroup scheme
$\sigma \subseteq G$ such that $X^{(\sigma)} \neq \emptyset$.
\end{definition}

A semicyclic subgroup scheme of $G$ is essential if and only if it is
essential for the action of $G$ on $X_s$ for some $s$; hence there are
only finitely many essential subgroups of $G$. We will fix a positive
integer $N$ that is divisible by the orders of all the essential subgroups
of $G$.

If $\sigma$ is a dual semicyclic subgroup of $G$, we define the
multiplicative system $\mathrm{S}_\sigma \subseteq \R(G/\sigma_0)\subseteq
\R G$ as before, as the subset of $\R (G/\sigma_0) \subseteq \R G$
consisting of elements whose image in $\Rtilde(\sigma/ \sigma_0)$ is a
power of $N$. Also,
    \[
    \K(X^{(\sigma)}, G)_\sigma  \eqdef \mathrm{S}_\sigma^{-1}\K(X^{(\sigma)},
    G),
    \]
as before.

\begin{proposition}\label{prop:inclusion}
Let $\sigma$ and $\tau$ be two semi-cyclic subgroups of $G$. If
$\tau\prec \sigma$, then\/ $\mathrm{S}_\sigma \subseteq
\mathrm{S}_\tau$.
\end{proposition}

\begin{proof}
Consider the commutative diagram of group schemes
    \[
    \xymatrix{
    G/\tau_0\ar@{->>}[rr]
    && G/\sigma_0 \\
    \tau/\tau_0\ar@{->>}[rr]\ar@{ >->}[u]
    && {}\sigma/\sigma_0 \; ;\ar@{ >->}[u]
    }
    \]
by taking representation rings we get a commutative diagram of rings
    \[
    \xymatrix{
    {}\R(G/\tau_0)\ar@{->>}[d]
    && {}\R(G/\sigma_0)\ar@{->>}[d]
    \ar@{ >->}[ll]\\
    {}\R(\tau/\tau_0)\ar@{->>}[d]
    && {}\R(\sigma/\sigma_0)\ar@{->>}[d]\ar@{ >->}[ll]
    \ar@{ >->}[ll]\\
    \Rtilde(\tau/\tau_0)
    && {}\Rtilde(\sigma/\sigma_0)\ar@{.>}[ll]
    }
    \]
(without the dotted arrow).
But it is easy to see that in fact the composition $\R(\sigma/\sigma_0)
\to \R(\tau/\tau_0) \to \Rtilde(\tau/\tau_0)$ factors through
$\Rtilde(\sigma/\sigma_0)$, so in fact the dotted arrows exists; and this
proves the thesis.
\end{proof}

Now, consider the restriction of the projection $\pi_{\sigma, \tau}\colon
\N_\sigma^{(\tau)} \to X^{(\sigma)}$, where $\sigma$ and $\tau$ are
dual semicyclic subgroups of $G$, and $\tau$ is subordinate to
$\sigma$. Because of Proposition~\ref{prop:inclusion}, we can consider the
composition of the pullback $\K(X^{(\sigma)}, G)_\sigma \to
\K(\N_\sigma^{(\tau)}, G)_\sigma$ with the natural homomorphism
$\K(\N_\sigma^{(\tau)}, G)_\sigma \to \K(\N_\sigma^{(\tau)}, G)_\tau$
coming from the inclusion $\mathrm{S}_\sigma \subseteq \mathrm{S}_\tau$; we
denote this homomorphism by
    \[
    \pi_{\sigma, \tau}^* \colon \K(X^{(\sigma)}, G)_\sigma
    \longrightarrow
    \K(\N_\sigma^{(\tau)}, G)_\tau.
    \]

\begin{definition}
Suppose that  $\sigma$ and $\tau$ are dual semicyclic subgroups of $G$ and
that $\tau$ is subordinate to $\sigma$. Two elements $a_\sigma \in
\K(X^{(\sigma)}, G)_\sigma$ and $a_\tau \in \K(X_\tau, G)_\tau$ are
\emph{compatible} if
    \[
    \pi_{\sigma, \tau}^*a_\sigma = \sp_{X, \sigma}^\tau
    a_\tau \in \K(\N_\sigma^{(\tau)}, G)_\tau.
    \]
\end{definition}

For each essential dual semicyclic subgroup $\sigma\subseteq G$ we denote
by
    \[
    \loc_\sigma \colon \K(X, G)\otimes\mathbb{Z}[1/N] \to \K(X^{(\sigma)},
    G)_\sigma
    \]
the composition of the restriction homomorphism
    \[\K(X, G) \otimes\mathbb{Z}[1/N] \to \K(X^{(\sigma)},G)
    \otimes\mathbb{Z}[1/N]
    \]
with the localization homomorphism
    \[
    \K(X^{(\sigma)},G) \otimes\mathbb{Z}[1/N] \to
    \K(X^{(\sigma)},G)_\sigma.
    \]

\bigskip

The following is the main result of this section.\\

\begin{theorem}\label{thm:refineddecomposition}
The ring homomorphism
    \[
    \prod_{\sigma}\loc_\sigma \colon \K(X, G) \otimes\mathbb{Z}[1/N]
    \longrightarrow\prod_{\substack{\sigma \subseteq G\\
    \sigma \text{\normalfont{} essential}}}
    \K(X^{(\sigma)}, G)_\sigma
    \]
is injective. Its image consists of the elements $(a_\sigma)$ of
$\prod_{\sigma} \K(X^{(\sigma)}, G)_\sigma$
with the property that if $\sigma$ and $\tau$ are essential,
$\tau\prec \sigma$, and $\dim \sigma = \dim \tau + 1$,
then $a_{\tau}$ and $a_\sigma$ are compatible.
\end{theorem}

Notice that in the particular case that the action has stabilizers of
constant dimension, all essential subgroups of $G$ have the same
dimension, and this reduces to Theorem~\ref{thm:refinedconstant}.

Also, if $\sigma$ is an essential subgroup of $G$ then
$X^{(\sigma)} = X_{\dim \sigma}^\sigma$, so it follows from Theorem
\refall{thm:refinedconstant}{2} that a splitting $G \simeq
(G/\sigma_0)
\times \sigma_0$ gives an isomorphism of rings
    \[
    \K(X^{(\sigma)}, G)_\sigma \simeq
    \K(X^{(\sigma)}, G/\sigma_0)\geom\otimes\Rtilde \sigma \otimes \R
    \sigma_0.
    \]
However, this isomorphism is not canonical in general, as it
depends on the choice of a splitting.

\begin{proof}
To simplify the notation, we will implicitly assume that everything has
been tensored with $\mathbb{Z}[1/N]$.

We apply Theorem~\ref{thm:maintheorem} together with
Theorem~\ref{thm:refinedconstant}. According to
Theorem~\ref{thm:maintheorem} we have an injection
$\K(X,G)\hookrightarrow \prod_{s}\K(X_{s},G)$ whose image is
the subring of sequences $(\alpha_s) \in \prod_{s=0}^n \K(X_s,G)$
with the property that for each $s= 1$, \dots,~$n$ the pullback of
$\alpha_s \in  \K(X_s,G)$ to $\K\bigl(\N_{s,s-1}, G\bigr)$ coincides
with\/ $\sp_{X,s}^{s-1}(\alpha_{s-1}) \in \K\bigl(\N_{s,s-1},
G\bigr)$. Moreover, by Theorem~\ref{thm:refinedconstant}, we can
decompose further each $\K(X_{s},G)$ as
$\prod_{\sigma}\K(X^{(\sigma}),G)_{\sigma}$, where $\sigma$
varies in the (finite) set of essential subgroups of $G$ of
dimension $s$. By compatibility of specializations, for
any $s\geq 0$, the following diagram is commutative
    \[
    \xymatrix{
    \K(X_{s-1}, G)\ar[r]^-\sim \ar[d]_{\prod_\tau\sp^{s-1}_{X,s}} &
    \prod\limits_{\substack{\tau \text{ essential}\\
    \dim \tau = s-1}}\K(X^{(\tau)}, G)_\tau
    \ar[d]^{\sp^\tau_{X,s}}\ar[dr]^{\prod_{\tau, \sigma}\sp^\tau_{X,
    \sigma}}\\
    \K(\N_{s,s-1}, G)\ar[r]^-\sim &
    \prod\limits_{\substack{\tau \text{ essential}\\
    \dim \tau = s-1}}\K(\N^{(\tau)}, G)_\tau\ar[r]^-\phi &
    \prod\limits_{\substack{\tau \text{ essential}\\
    \dim \tau = s-1}}\prod\limits_{\substack{\sigma\succ \tau\\
    \dim \sigma = s}}\K(\N_\sigma^{(\tau)}, G)
    }
    \]
where $\phi$ is induced by the obvious pullbacks. On the other
hand, the following diagram commutes by definition of
$\pi_{\sigma, \tau}^*$
\[
\xymatrix{
\K(X_s, G)\ar[r]^-\sim \ar[d]&
\prod\limits_{\substack{\sigma \text{ essential}\\
\dim \sigma = s}}\K(X^{(\sigma)}, G)_\sigma
\ar[dr]^{\prod_{\sigma, \tau}\pi^*_{\sigma, \tau}}\\
\K(\N_{s,s-1}, G)\ar[r]^-\sim &
\prod\limits_{\substack{\tau \text{ essential}\\
\dim \tau = s-1}}\K(\N^{(\tau)}, G)_\tau\ar[r]^-\phi &
\prod\limits_{\substack{\tau \text{ essential}\\
\dim \tau = s-1}}\prod\limits_{\substack{\sigma\succ \tau\\
\dim \sigma = s}}\K(\N_\sigma^{(\tau)}, G)
}.
\]
Then the Theorem will immediately follow if we show
that $\phi$ is an isomorphism. This is true because of the following
Lemma.

\begin{lemma}
  Fix an essential subgroup $\tau$ of dimension $s-1$. Then for any $\sigma
\succ \tau$ with $\dim \sigma = s$, the scheme $\N_\sigma^{(\tau)}$ is open
in $\N_{s}^{(\tau)}$; furthermore, $\N_{s}^{(\tau)}$ is the disjoint union of
the $\N_\sigma^{(\tau)}$ for all essential $\sigma$ with $\sigma \succ
\tau$, $\dim \sigma = s$.
\end{lemma}

\begin{proof}
We will show that $\N_{s}^{(\tau)}$ is the disjoint union of the
$\N_\sigma^{(\tau)}$; since each $\N_\sigma^{(\tau)}$ is closed in
$\N_s^{(\tau)}$, and there are only finitely many possible $\sigma$, it
follows that each $\N_\sigma^{(\tau)}$ is also open in $\N_{s}^{(\tau)}$.

Let us first observe that if $\sigma$ and $\sigma'$
are essential subgroups in $G$ of dimension $s$ to which $\tau$ is
subordinate, and
$\N_{\sigma}^{(\tau)}\bigcap\N_{\sigma'}^{(\tau)} \neq \emptyset$, then
$X^{(\sigma)}\bigcap X^{(\sigma')}\neq \emptyset$, therefore
$\sigma_{0}$ is equal to $\sigma'_{0}$. But this implies that
$\sigma=\sigma'$, since $\sigma$ and $\sigma'$ are both equal to the
inverse image in
$G$ of the image of $\tau\rightarrow G{/}\sigma_{0} = G{/}\sigma'_{0}$.

According to Proposition~\refall{prop:firststratification}{3}, if $T_1$,
\dots,~$T_r$ are the essential $s$-dimensional subtori of $G$, then $\N_s$
is the disjoint union of the $\N_{T_j}$. Clearly, if $\tau_0$ is not
contained in $T_j$, then $\N_{T_j}^{(\tau)}$ is empty, so $\N_s^{(\tau)}$
is the disjoint union of the $\N_T^{(\tau)}$, where $T$ ranges over the
essential $s$-dimensional subtori of $G$ with $\tau_0 \subseteq T$. But
there is a bijective correspondence between $s$-dimensional dual
semi-cyclic subgroups $\sigma \subseteq G$ with $\sigma \succ \tau$ and
$s$-dimensional subtori of $G$ with $\tau_0 \subseteq T$: in one direction
we associate with each $\sigma$ its toral part $\sigma_0$, in the other we
associate with each $T$ the subgroup scheme $\tau + T \subseteq G$.

The proof is concluded by noticing that if $\sigma$ and $\tau$ are as
above, with $\sigma_0 = T$, then $\N_T^{(\tau)} = \N_\sigma^{(\tau)}$.
\end{proof}
\noqed\end{proof}

\section*{Addendum: corrections (August 2004)}

Amnon Neeman has noticed a serious error in the proof of Theorem~ 3.2: the  argument given does not yield a uniquely defined specialization map $\mathrm{Sp}_{Y}$, so that in particular compatibility with pullbacks does not hold. This is due to the elementary fact, overlooked in the paper, that if one has a fiber sequence of spectra $$\xymatrix{E \ar[r]^-{f} & E' \ar[r]^-{g} & E'' },$$ then a map $h:E'\arr E$ such that $h\circ f$ is homotopic to zero, does induce a map of spectra $p:E''\arr E$, but this map is not unique, as it can be modified by using any map $E[1]\rightarrow E$ (by adding to any given $p$ the composite $E''\rightarrow E[1] \rightarrow E$). Of course if the nullhomotopy  $h\circ f \sim 0$ is specified then this singles out a unique map $p:E''\rightarrow E$; but it is not clear to the authors how to choose a homotopy; thus, they are unable to define a specialization map in the  generality claimed in the statement of Theorem~3.2.

Fortunately, it is still possible to define the specialization homomorphisms for higher K-theory in a generality that is sufficient for the rest of the paper: thus, all the results in sections 4, 5, 6 and 7, including the two main theorems, still hold unchanged. Also, section 2 which is independent of section 3 where specializations were defined,  remains unchanged.\\

In what follows we will work in the same setup as in the paper, to which we refer for the unexplained notation. If $Y$ is a closed subscheme of a scheme $X$ over a fixed base $S$, we denote by $\rM^{0}_{Y}X \arr \PP^{1}_{S}$ the deformation to the normal bundle, as in \cite[Chapter~5]{fulton}, and in the paper. We denote by $\infty$ the closed subscheme of $\PP^{1}_{S}$ that is the image of the section at infinity $S \arr \PP^{1}_{S}$; the inverse image of $\infty$ in $\rM^{0}_{Y}X$ is the normal bundle $\N_{Y}X$.

Assume that $X$ is a regular noetherian algebraic space with the action of a diagonalizable group $G$, $Z$ a $G$-invariant regular Cartier divisor with trivial normal bundle, $i \colon Z \into X$  and $j \colon X \setminus Z \subseteq X$ the embeddings. The composition
   \[
   \K(Z,G) \stackrel{i_{*}} \longrightarrow
      \K(X, G) \stackrel{i^{*}} \longrightarrow \K(Z,G)
   \]
is $0$: if we assume that $j^{*} \colon \K(X,G) \arr \K(X \setminus Z, G)$ is surjective, then we have an exact sequence
   \[
   0 \arr \K(Z,G) \stackrel{i_{*}} \longrightarrow
      \K(X, G) \stackrel{j^{*}} \longrightarrow 
      \K(X \setminus Z, G) \arr 0;
   \]
hence the homomorphism $i^{*}\colon \K(X,G) \arr \K(Z,G)$ factors through $\K(X \setminus Z, G)$, inducing a specialization ring homomorphism
   \[
   \sp^{X}_{Z}\colon \K(X \setminus Z, G) \arr \K(Z,G).
   \]
If we restrict to $\K[0]$, then surjectivity holds, and this is already in \cite[X-Appendice, 7.10]{sga6}.

Recall that $X_{s}$ is the regular subscheme of $X$ where the stabilizers have fixed dimension $s$, and that we have set $\rM_s \eqdef \rM^{0}_{X_{s}}X \arr \PP^{1}$. Consider the closed embedding $\N_{s} \subseteq \rM_{s}$, whose complement is $X_{\leq s} \times \AA^{1}$. Looking at the composition
   \[
   X_{\leq s} \times \AA^{1} \into \rM_{s} \arr X_{\leq s} \times \PP^{1} \arr X_{\leq s}
   \]
we see that the pullback $\K(\rM_{s}, G) \arr \K(X_{\leq s}, G)$ is surjective. Consider now the open embedding $X_{\leq t} \subseteq X_{\leq s}$ (for $s\geq t$): the pullback $\K(X_{\leq s}, G) \arr \K(X_{\leq t}, G)$ is also surjective, by K-rigidity. From the commutative diagram
   \[
   \xymatrix{
   \K(\rM_{s}, G) \ar@{->>}[r]\ar[d] & \K(X_{\leq s}, G)\ar@{->>}[d]\\
   \K(\rM_{s, \leq t}, G) \ar[r] & \K(X_{\leq t}, G)
   }
   \]
we conclude that the restriction $\K(\rM_{s, \leq t}, G) \arr \K(X_{\leq t}, G)$ is surjective, so we have an exact sequence
   \[
   0 \arr \K\bigl(\N_{s, \leq t}, G\bigr) \arr \K(\rM_{s, \leq t}, G) \arr \K(X_{\leq t}, G) \arr 0.
   \]
This allows to define specialization maps
   \[
   \sp^{\leq t}_{X,s}\eqdef\sp^{\rM_{s, \leq t}}_{\N_{s, \leq t}} \colon \K(X_{\leq t}, G) \arr \K\bigl(\N_{s, \leq t}, G\bigr).
   \]

To define $\sp^{t}_{X,s}$ consider the commutative diagram with exact rows
   \[
   \xymatrix{
   0\ar[r] & \K(X_{t},G)\ar[r]\ar@{-->}[d] &  
      \K(X_{\leq t}, G)\ar[r]\ar[d]^{\sp^{\leq t}_{X,s}} &
      \K(X_{\leq t-1}), G\ar[d]^{\sp^{\leq t-1}_{X,s}}\ar[r] & 0 \\
   0\ar[r] & \K\bigl(\N_{s, t},G \bigr)\ar[r] &
      \K\bigl(\N_{s, \leq t},G \bigr) \ar[r] &
      \K\bigl(\N_{s, \leq t-1},G \bigr)\ar[r] & 0
   }
   \]
(the commutativity of the second square follows easily from functoriality of pullbacks).

\begin{definition}\label{def:specialization1}
The specialization homomorphism
   \[
   \sp^{t}_{X,s} \colon  \K(X_{t},G) \arr \K\bigl(\N_{s, t},G \bigr)
   \]
is the unique dotted arrow that fits in the diagram above.
\end{definition}

These coincide with the usual specialization homomorphisms for $\K[0]$. This is clear for the $\sp^{\leq t}_{X.s}$. For $\sp^{t}_{X,s}$ it follows from the fact the cartesian diagram
   \[
   \xymatrix{
   \N_{s,t} \ar[r]\ar[d] & \rM_{s,t} \ar[d]\\
   \N_{s,\leq t} \ar[r]s & \rM_{s,\leq t}
   }
   \]
is Tor-independent, and from the following Lemma.

\begin{lemma}\label{lem:tor-independent}
If 
   \[
   \xymatrix{
   X' \ar[r]^{f'} \ar[d]^{\phi} & Y' \ar[d]^{\psi}\\
   X  \ar[r]^{f}                & Y
   }
   \]
is a Tor-independent cartesian square of regular algebraic spaces with an action of $G$, where $f$ is a closed embedding. Then the diagram
   \[
   \xymatrix{
   \K(X,G) \ar[r]^{f_{*}} \ar[d]^{\phi^{*}} & \K(Y,G) \ar[d]^{\psi^{*}}\\
   \K(X',G)  \ar[r]^{f'_{*}}                & \K(Y',G)
   }
   \]
commutes.
\end{lemma}

The proof that starts at the top of page~10 is general enough.

Now we have to check compatibility of specializations.

\begin{lemma}
Denote by $i$ the inclusion of $X_{t}$ in $X_{\leq t}$ and by $i'$ that of $\N_{s,t}$ in $\N_{s,\leq t}$. Then the diagram
   \[
   \xymatrix{
   \K(X_{\leq t}, G) \ar[r]^-{i^{*}} \ar[d]^{\sp^{\leq t}_{X, s}} &
      \K(X_{t}, G) \ar[d]^{\sp^{t}_{X, s}}\\
   \K(\N_{s,\leq t}, G) \ar[r]^-{i'^{*}} &
      \K(\N_{s,t}, G)
   }
   \]
commutes.
\end{lemma}

\begin{proof}
By the definition of $\sp^{t}_{X,s}$, we need to check that
   \[
   \xymatrix{
   \K(X_{\leq t}, G) \ar[r]^-{i^{*}} \ar[d]^{\sp^{\leq t}_{X, s}} &
      \K(X_{t}, G) \ar[r]^-{i_{*}} \ar[d]^{\sp^{t}_{X, s}} &
      \K(X_{\leq t}, G) \ar[d]^{\sp^{\leq t}_{X, s}}\\
   \K(\N_{s,\leq t}, G) \ar[r]^-{i'^{*}} &
      \K(\N_{s,t}, G) \ar[r]^-{i'_{*}}&
      \K(\N_{s,\leq t}, G)
   }
   \]
commutes. By the projection formula (see \cite[Proposition~A.5]{vevi}) we see that the group homomorphisms $i_{*}i^{*}$ and $i'_{*}i'^{*}$ are multiplications by 
   \[
   [i_{*}\cO_{X_{t}}] \in \K((X_{\leq t}, G) \quad\text{and}\quad
   [i_{*}\cO_{\N_{s,t}}] \in \K(\N_{s,\leq t}, G)
   \]
respectively: so we have to prove that the diagram
   \[
   \xymatrix@C+15pt{
   \K(X_{\leq t}, G) \ar[r]^-{\cdot [i_{*}\cO_{X_{t}}]}
                     \ar[d]^{\sp^{\leq t}_{X, s}} &
      \K(X_{\leq t}, G) \ar[d]^{\sp^{\leq t}_{X, s}}\\
   \K(\N_{s,\leq t}, G) \ar[r]^-{\cdot [i_{*}\cO_{\N_{s,t}}]} &
      \K(\N_{s,\leq t}, G)
   }
   \]
commutes. Since $\sp^{\leq t}_{X, s}$ is a ring homomorphism, this is equivalent to saying that
   \[
   \sp^{\leq t}_{X, s}[i_{*}\cO_{X_{t}}] = 
      [i_{*}\cO_{\N_{s,t}}] \in \K[0](\N_{s,\leq t}, G).
   \]
But $[i_{*}\cO_{X_{t}}]$ is the restriction of $[i_{*}\cO_{\rM_{s,t}}] \in \K[0](\rM_{s,\leq t}, G)$, so we have to show that the restriction of $[i_{*}\cO_{\rM_{s,t}}]$ to $\N_{s, \leq t}$ is $[i_{*}\cO_{\N_{s,t}}]$; and this follows immediately from the fact that the square
   \[
   \xymatrix{
   \N_{s,t} \ar[r]\ar[d] & \rM_{s,t} \ar[d]\\
   \N_{s,\leq t} \ar[r]  & \rM_{s,\leq t}
   }
   \]
is cartesian and Tor-independent, by Lemma~\ref{lem:tor-independent}.
\end{proof}

With this definition, and the compatibility property proved above, everything goes through in Sections 4, 5 and 6. For the theory of Section~7 to work, we need to define specialization maps
   \[
   \K(X^{(\tau)},G) \arr \K(\N_{s}^{(\tau)}, G)
   \]
when $\tau$ is a diagonalizable subgroup scheme of $G$ and $s$ is an integer with $s \geq \dim\tau$ (see the bottom of page~39).The cartesian diagram of embeddings
   \[
   \xymatrix{
   X_s^{\tau} \ar[r]\ar[d] & X^{\tau} \ar[d] \\
   X_s        \ar[r]       & X \\
   }
   \]
yields an embedding of $G$-spaces $\N_{X_s^{\tau}}X^{\tau} \into \N_s$; since $\tau$ acts trivially on $\N_{X_s^{\tau}}X^{\tau}$ we get an embedding $\N_{X_s^{\tau}}X^{\tau} \into \N_s^{\tau}$.

\begin{lemma}\label{lem:emb-isom}
The embedding $\N_{X_s^{\tau}}X^{\tau} \into \N_s^{\tau}$ is an isomorphism.
\end{lemma}

\begin{proof}
Consider the natural embedding of deformations to the normal bundle
   \[
   \xymatrix{
   \rM^{0}_{X_{s}^{\tau}}X^{\tau} \into 
   (\rM_{s})^{\tau};
   }
   \]
generically, that is, over $\AA^{1}$, they coincide. On the other hand it follows from Proposition~3.6 in the paper that the inverse image of $\AA^{1}$ in $(\rM_{s})^{\tau}$ is scheme-theoretically dense in $(\rM_{s})^{\tau}$, and this shows that this embedding is an isomorphism. Since the fibers over $\infty$ of $\rM^{0}_{X_{s}^{\tau}}X^{\tau}$ and $(\rM_{s})^{\tau}$ are $\N_{X_s^{\tau}}X^{\tau}$ and $\N_s^{\tau}$ respectively, this concludes the proof.
\end{proof}

Now set $t = \dim \tau$, so that $X^{(\tau)} \eqdef X^{\tau}_{t}$. Then we get a specialization map
   \[
   \sp_{X,s}^{\tau} \eqdef \sp_{X^{\tau},s}^{t} \colon \K(X^{\tau}_{t}, G) = \K(X^{(\tau)}, G) \arr
                         \K\bigl((\N_{s}^{\tau})_{t}, G\bigr) = \K(\N_{s}^{(\tau)}, G)
   \]
that is exactly what we want. This allows to define the specialization map
   \[
   \sp_{X, \sigma}^{\tau} \colon \K(X^{\tau}_{t}, G) \arr
      \K(\N_{\sigma}^{(\tau)}, G)
   \]
for any pair of dual cyclic subgroups $\sigma$ and $\tau$, with $\tau \prec \sigma$, as on page~41 of the paper, by composing $\sp_{X^{\tau},s}^{t}$ with the restriction homomorphism $\K(\N_{s}^{(\tau)}, G) \arr \K(\N_{\sigma}^{(\tau)}, G)$. Note that $X^{(\tau)}=X^{\tau}_{\leq t}$, $\N_{s}^{(\tau)}=(\N^{\tau}_{s})_{\leq t}$,  and $\sp_{X,s}^{\tau}$ can also be identified with $\sp_{X^{\tau},s}^{\leq t}$; therefore $\sp_{X,s}^{\tau}$ and $\sp_{X, \sigma}^{\tau}$ are ring homomorphisms.

\subsection*{Further corrections} Here we correct a few typos that we have noticed since the publication of the article.

In  the statement of Proposition~1.1, ``scheme'' should be replaced by ``algebraic space''.

The are several typos in the diagrams on p. 42:

\begin{enumerate}

\item $\prod \sp ^{s-1}_{X,s}$ should be replaced by $\sp ^{s-1}_{X,s}$,

\item $\sp ^{\tau}_{X,s}$ by $\prod \sp ^{\tau}_{X,s}$,

\item $\N^{(\tau)}$ by $\N^{(\tau)}_{s}$ and

\item $\K(\N^{(\tau)}_{\sigma},G)$ by $\K(\N^{(\tau)}_{\sigma},G)_{\tau}$

\end{enumerate}

Finally, in the statement of Lemma~4.9, ``linearly independent elements'' should read ``pairwise linearly independent elements'' (we owe this also to Amnon Neeman).

\subsection*{Acknowledgments} We are very much in debt with Amnon Neeman who read our paper carefully and kindly pointed out the problem to us.


\begin{thebibliography}{B-B-F-K02}

\bibitem [Ati74]{atiyah} M.F. Atiyah: \emph{Elliptic
operators and compact groups},  Lecture Notes in Mathematics \textbf{401},
Springer--Verlag, Berlin--New York (1974).

\bibitem[B-B-F-K02]{bbfk} G. Barthel, J.P. Brasselet, K.H. Fieseler, L. Kaup:
\emph{Combinatorial intersection cohomology for fans}, Tohoku Math. J.
\textbf{54} (2002), 1--41.

\bibitem[Bia73]{b-b} A. Bia\l ynicki-Birula, \emph{Some theorems on
actions of algebraic groups}, Ann. of Math. \textbf{98} (1973),
480--497.

\bibitem[B-DC-P90]{bdcp} E. Bifet, C. De Concini, C. Procesi:
\emph{Cohomology of regular embeddings}, Adv. Math. \textbf{82} (1990),
1--34.

\bibitem[Bre74]{bredon} G.E. Bredon,
\emph{The free part of a torus action and related numerical equalities},
Duke Math. J. \textbf{41} (1974), 843--854.

\bibitem[Bri97]{brion1} M. Brion, \emph{Equivariant Chow groups for torus
actions}, Transform. Groups 2 (1997), 225--267.

\bibitem[Bri98]{brion2} \bysame, \emph{Equivariant
cohomology and equivariant intersection theory}, in \emph{Representation
theories and algebraic geometry (Montreal, PQ, 1997)}, Kluwer Acad.
Publ., Dordrecht, (1998), 1--37, .

\bibitem[Ch-Sk74]{ch-sk}T. Chang, T. Skjelbred, \emph{The topological
Shur lemma and related results}, Annals of Mathematics \textbf{100} (1974),
307--321.

\bibitem[Dan78]{danilov} V. Danilov, \emph{The geometry of toric
varieties}, Russ. Math. Surveys \textbf{33} (1978), 97--154.

\bibitem[SGA3]{sga3} M. Demazure, \emph{Sch\'emas en groupes},
Lecture Notes in Mathematics \textbf{151}, \textbf{152} and \textbf{153},
Springer-Verlag (1970).

\bibitem[Ed-Gr98]{edidingraham} D. Edidin, W. Graham, \emph{Equivariant
intersection theory}, Invent. Math. \textbf{131} (1998),
595--634.

\bibitem[Ei-Mo62]{em} S. Eilenberg, J.C. Moore: \emph{Limits and
spectral sequences}, Topology \textbf{1} (1962), 1--23.

\bibitem[Ful93a]{fulton} W. Fulton, \emph{Intersection Theory},
Springer-Verlag, (1993).

\bibitem[Ful93b]{fultontoric} W. Fulton, \emph{Introduction to toric
varieties}, Princeton University Press (1993).

\bibitem[G-K-MP98]{gkmp} M. Goresky, R. Kottwitz, R. MacPherson,
\emph{Equivariant cohomology, Koszul duality, and the localization
theorem}, Invent. Math. \textbf{131} (1998), 25--84.

\bibitem[EGAIV]{ega4} A. Grothendieck, J. Dieudonne, \emph{Elements de
G\`eom\`etrie Alg\`ebrique IV. . \'Etude locale des sch\'emas et des
morphismes de sch\'emas}, Publ. Math. IHES \textbf{28} (1966).

\bibitem[Hsi75]{hsiang} W.\,Y. Hsiang, \emph{Cohomology theory of
Topological Transformation Groups}, Springer Verlag (1975).

\bibitem[Kir84]{kirwan} F. Kirwan, \emph{Cohomology of quotients in
symplectic and algebraic geometry}, Mathematical Notes, \textbf{31},
Princeton University Press, Princeton (1984).

\bibitem[Kly83]{kly} A.A. Klyachko, \emph{Vector bundles on Demazure
models}, Selecta Math. Soviet. \textbf{3} (1983/84), 41--44.

\bibitem[La-MB00]{lmb} G. Laumon, L. Moret-Bailly, \emph{Champs
Alg\'ebriques}, Springer--Verlag (2000).

\bibitem[Mer97]{merk} A. Merkurjev, \emph{Comparison of the
equivariant and the standard K-theory of algebraic varieties.} Algebra i
Analiz \textbf{9} (1997),175--214. Translation in St. Petersburg
Math. J. \textbf{9} (1998), 815--850.

\bibitem[Ro-Kn]{rokn} I. Rosu, A. Knutson, Equivariant K-theory and
Equivariant Cohomology, preprint	math.AT/9912088, to appear in
Math. Zeit.

\bibitem[SGA6]{sga6} P. Berthelot, A. Grothendieck, L. Illusie, \textit{Th\'eorie des intersections et th\'eor\`eme de Riemann-Roch}, LNM 225, Springer, Berlin 1971.

\bibitem[Sum75]{sumihiro} H. Sumihiro, \emph{Equivariant completion
II}, J. Math. Kyoto Univ. \textbf{15} (1975), 573--605.

\bibitem[Tho86a]{th2} R.\,W. Thomason, \emph{Comparison of equivariant
algebraic and topological K-theory}, Duke Math. J. \textbf{68}
(1986), 447--462.

\bibitem[Tho86b]{th4}\bysame, \emph{Lefschetz--Riemann--Roch theorem and
coherent trace formula}, Invent. Math. \textbf{85} (1986),
515--543.

\bibitem[Tho87]{th5} \bysame, \emph{Algebraic K-theory of group scheme
actions}, in \emph{Algebraic Topology and Algebraic K-theory}, Annals of
Mathematicals Studies \textbf{113}, Princeton University Press, Princeton
(1987)

\bibitem[Tho92]{th3} \bysame, \emph{Un formule de Lefschetz en
K-th\`eorie equivariante alg\`ebrique}, Duke Math. J. \textbf{68}
(1992), 447--462.

\bibitem[Tho93]{th1} \bysame, \emph{Les K-groupes d'un sch\'ema
\'eclat\'e et un formule d'intersection exc\'eden\-taire}, Invent. Math.
\textbf{112} (1993), 195--216.

\bibitem[Th-Tr90]{thtr} \bysame, T. Trobaugh, \emph{Higher
algebraic K-theory of schemes and of derived categories}, Grothendieck
Festschrift vol. III, Birkh\"auser (1990), 247--435.

\bibitem[Toen]{toen} B. Toen, \emph{Notes on \textrm{G}-theory of
Deligne--Mumford stacks}, preprint math.AG/9912172.

\bibitem[Ve-Vi]{vevi} G. Vezzosi, A. Vistoli, \emph{Higher algebraic
K-theory of group actions with finite stabilizers}, Duke Math. J., \textbf{113}
(2002), 1--55.




\end{thebibliography}
\end{document}